\documentclass[11pt,reqno]{amsart}

\usepackage{geometry}
\usepackage{graphicx,float}
\usepackage[hidelinks]{hyperref}
 \usepackage{booktabs} 

\RequirePackage{amsthm,amsmath,amsfonts,amssymb,dsfont,xcolor}

\newif\ifshowchanges
\showchangestrue

\theoremstyle{plain}
\newtheorem{theorem}{Theorem}[section]
\newtheorem{proposition}{Proposition}[section]
\newtheorem{lemma}{Lemma}[section]

\theoremstyle{definition}

\newtheorem{remark}{Remark}[section]

\newcommand{\mathbbm}[1]{\mathds{#1}}
\newcommand{\R}{\mathbb{R}}

\newcommand{\N}{\mathbb{N}}
\newcommand{\E}{\mathbb{E}}
\newcommand{\Prob}{\mathbb{P}}

\newcommand{\calB}{\mathcal{B}}

\newcommand{\calX}{\mathcal{X}}
\newcommand{\calY}{\mathcal{Y}}
\newcommand{\calC}{\mathcal{C}}
\newcommand{\calP}{\mathcal{P}}
\newcommand{\calR}{\mathcal{R}}
\newcommand{\calD}{\mathcal{D}}
\newcommand{\calS}{\mathcal{S}}

\newcommand{\calQ}{\mathcal{Q}}

\newcommand{\sfA}{\mathsf{A}}
\newcommand{\sfX}{\mathsf{X}}
\newcommand{\sfY}{\mathsf{Y}}
\newcommand{\sfc}{\mathsf{c}}
\newcommand{\sfC}{\mathsf{C}}

\newcommand{\frakq}{\mathfrak{q}}
\newcommand{\fraku}{\mathfrak{u}}

\newcommand{\GW}{\mathsf{GW}_{2,2}}
\newcommand{\Wp}{\mathsf{W}_p}
\newcommand{\ind}{\mathbbm{1}}

\newcommand{\TV}{\mathsf{d}_{\mathrm{TV}}}

\newcommand{\Dminus}{\calD_{M,\eta}^-}
\newcommand{\Dplus}{\calD_{M}^+}

\DeclareMathOperator{\dist}{dist}

\DeclareMathOperator{\Var}{Var}
\DeclareMathOperator{\Cov}{Cov}

\DeclareMathOperator{\spt}{spt}

\begin{document}
\title[Minimax rates for OT distances]{Minimax rates  for estimating optimal transport distances with vanishing separation}
\thanks{K. Kato is partially supported by NSF grant DMS-2413405.}
\date{First version: November 19, 2026. This version: \today.}

\author[K. Kato]{Kengo Kato}
\address[K. Kato]{
Department of Statistics and Data Science, Cornell University.
}
\email{kk976@cornell.edu}

\author[B. Wang]{Boyu Wang}
\address[B. Wang]{
Department of Statistics and Data Science, Cornell University.
}
\email{bw563@cornell.edu}

\begin{abstract}
We study minimax rates of convergence for estimating the $p$-Wasserstein distance $\Wp$ for $p \ge 1$ and the quadratic Gromov--Wasserstein distance $\GW$ under a separation constraint where the population distance is at least $\underline{\delta} > 0$. For $\Wp$, we consider compactly supported marginals on a common Euclidean space, $\R^d$, and for $\GW$, we consider compactly supported marginals on Euclidean spaces with possibly different dimensions. Previous work has shown that when $\underline{\delta}$ is bounded away from zero, the empirical estimator, based on $n$ independent observations from each of the marginals, may enjoy faster convergence rates than under the unrestricted case. In this paper, we establish sharp joint dependence of the minimax rate on $(n,\underline{\delta})$, up to logarithmic factors, where $\underline{\delta} = \underline{\delta}_n \downarrow 0$ is allowed. For $\Wp$, our results reveal new phase transition phenomena in the dependence on $\underline{\delta}$.
For both problems, the empirical estimator attains the upper bounds. For several configurations of $(p,d)$ under $\Wp$ (e.g., $d > 2p > 4$), we establish improved error bounds in $\underline{\delta}$ for the empirical estimator. The improved dependence follows from novel error bounds for the empirical transport cost that involve the size of the population distance, which may be of independent interest. 

\end{abstract}

\maketitle

\section{Introduction}
\subsection{Overview}
Estimating optimal transport (OT) distances constitutes a fundamental problem in statistics, machine learning, and probability \cite{villani2009optimal,santambrogio2015optimal}. 
We refer to \cite{chewi2025statistical} for an overview of statistical OT. A feature of distance estimation is that its difficulty can depend on the size of the distance being estimated. For Wasserstein distances, the empirical estimator can converge faster when the population distance is bounded away from zero than the rate available uniformly over all pairs of distributions. Understanding the transition between these regimes requires bounds that track the population distance as it approaches zero. In this work, we study minimax rates of convergence for estimating the $p$-Wasserstein distance $\Wp$ for $p \ge 1$ and the quadratic Gromov--Wasserstein (GW) distance $\GW$ \cite{memoli2011gromov,sturm2023the}, under a separation constraint where the distance is at least $\underline{\delta}$. Here, we allow $\underline{\delta} \downarrow 0$ as the sample size increases. 

To motivate our results, consider estimating $\Wp(\mu,\nu)$ from $n$ independent observations from each of compactly supported Borel probability measures $\mu,\nu$ on $\R^d$ (see Section \ref{sec: preliminaries} for the definitions of $\Wp$ and $\GW$ and the formal setting). A natural estimator for $\Wp$ is the empirical estimator, $\Wp(\hat{\mu}_n,\hat{\nu}_n)$, where $\hat{\mu}_n,\hat{\nu}_n$ are the empirical distributions for $\mu,\nu$. For now, let $\lesssim$ denote an inequality that holds up to a constant independent of $(n,\underline{\delta})$.  Theorem 1 in \cite{fournier2015on} directly implies, for $d > 2p$, 
\begin{equation}
\E[|\Wp(\hat{\mu}_n,\hat{\nu}_n) - \Wp(\mu,\nu)|] \lesssim n^{-1/d},
\label{eq: baseline}
\end{equation}
which agrees with the minimax lower bound from Theorem 11 in \cite{nilesweed2022estimation}, up to logarithmic factors, for the class of distributions supported in a fixed ball. \cite{chizat2020faster,manole2024sharp} made an interesting observation that faster rates are possible when $\Wp(\mu,\nu)$ is bounded away from zero. Indeed, \cite{chizat2020faster,manole2024sharp} established that, for $d > 2p$,
\begin{equation*}
\E[|\Wp^p(\hat{\mu}_n,\hat{\nu}_n)-\Wp^p(\mu,\nu)|] \lesssim n^{-(p \wedge 2)/d},
\label{eq: known transport rate}
\end{equation*}
which directly implies that, whenever $\Wp(\mu,\nu) \ge \underline{\delta} > 0$,
\begin{equation}
\E[|\Wp(\hat{\mu}_n,\hat{\nu}_n)-\Wp(\mu,\nu)|] \lesssim \underline{\delta}^{1-p} n^{-(p \wedge 2)/d}.
\label{eq: known rate}
\end{equation}
See Corollary 4 and Remark 1 in \cite{manole2024sharp}; see also \cite{hundrieser2024empirical}.
 One can take the minimum of two rates to obtain the rate $\min\{n^{-1/d}, \underline{\delta}^{1-p}n^{-(p \wedge 2)/d}\}$ for $\Wp(\hat{\mu}_n,\hat{\nu}_n)$. In particular, for $p \ge 2$, when $\underline{\delta}>0$ stays fixed, the rate improves to $n^{-2/d}$ from $n^{-1/d}$. A natural question, which has not been fully addressed in the literature, would be:
\begin{quote}
Does the rate $\min\{n^{-1/d}, \underline{\delta}^{1-p}n^{-(p \wedge 2)/d}\}$ match a minimax lower bound, possibly up to logarithmic factors, for a sufficiently large class of distributions $\mu,\nu$ satisfying $\Wp(\mu,\nu) \ge \underline{\delta} > 0$?
\end{quote}
The goal of the present paper is to explore sharp joint dependence of the minimax rate on $(n,\underline{\delta})$ for $\Wp$, including the $d \le 2p$ case, and a relevant problem for $\GW$. 

\begin{figure}[H]
    \centering
    \includegraphics[width=0.7\linewidth]{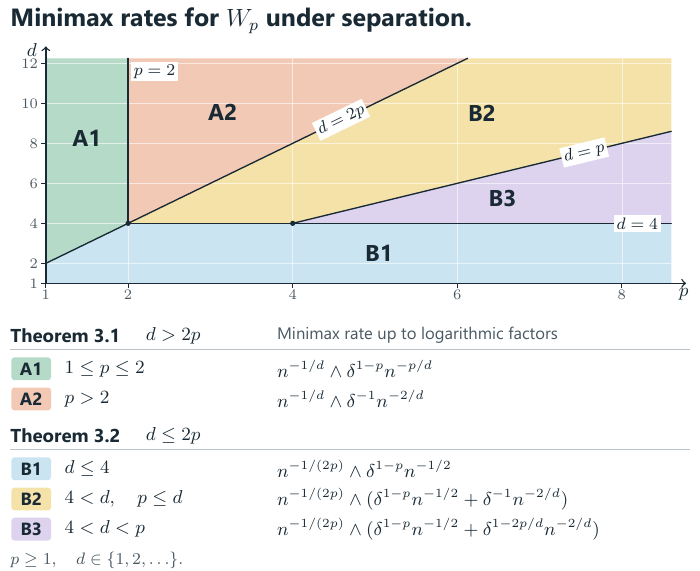}
    \caption{Minimax rates for $\Wp$ in Theorems~\ref{thm: Wp main} and~\ref{thm: Wp main small p}, up to logarithmic factors. The figure uses $\delta$ for the separation parameter $\underline{\delta}$.}
    \label{fig: wp-phase-regions}
\end{figure}

Our results for $\Wp$ are summarized in Figure \ref{fig: wp-phase-regions}. For all cases, the empirical estimator, without requiring the knowledge of $\underline{\delta}$, attains the upper bounds.  For $d > 2p$, the minimax rate is
\[
n^{-1/d}\wedge\bigl(\delta^{1-(p \wedge 2)}n^{-(p \wedge 2)/d}\bigr),
\]
up to logarithmic factors, uniformly over sufficiently small $\underline{\delta} > 0$. For $1\le p\le2$, the result establishes optimality, up to logarithmic factors, of the preceding empirical upper bound. At $p = 2$, the minimax rate exhibits a phase transition in $\underline{\delta}$, from $\underline{\delta}^{1-p} \ (1 \le p \le 2)$ to $\underline{\delta}^{-1} \ (p > 2$). For $p > 2$, we establish that the empirical estimator attains the improved rate in $\underline{\delta}$. To this end, we prove a novel error bound for the empirical \textit{transport cost}, $\Wp^p(\hat{\mu}_n,\hat{\nu}_n)$: for $t = \Wp(\mu,\nu)$ and $d > 2p > 4$, 
\begin{equation}
\E[|\Wp^p(\hat{\mu}_n,\hat{\nu}_n) -t^p|] \lesssim t^{p-2}n^{-2/d} + n^{-p/d},
\label{eq: transport-cost bound}
\end{equation}
which involves the size of the population distance (the inequality $\lesssim$ here does not depend on $t$). The bound implies that small $\Wp$ may lead to faster rates for $\Wp^p(\hat{\mu}_n,\hat{\nu}_n)$, which may be of independent interest. 

For $d \le 2p$, a global picture is more complicated. One striking result, at least from the authors' perspective, is a phase transition phenomenon between $p \le d$ and $p > d$, under $d > 4$. Under $d > 4$, the minimax rate involves $\delta^{-1} n^{-2/d}$ when $p \le d$, while it involves $\delta^{1-2p/d}n^{-2/d}$ when $p > d$.
Such a phase transition occurring at $d = p$ is new in the OT literature, to the best of the authors' knowledge. 
Proving sharp rates for the empirical estimator for $4 < d < p$, both in $(n,\underline{\delta})$, requires extra effort and rests on establishing $O(n^{-1/2} + t^{p(d-2)/d}n^{-2/d}\log n)$ rates for the empirical transport cost, $\Wp^p(\hat{\mu}_n,\hat{\nu}_n)$. 

For the quadratic GW distance, the marginals may be defined on spaces with different dimensions, $\R^{d_x}$ and $\R^{d_y}$. Theorem \ref{thm: GW main} shows that the minimax rate for $\GW$ under $\GW \ge \underline{\delta}$ is, up to logarithmic factors,
\[
n^{-1/(d \vee 4)} \wedge \big(\delta^{-1} n^{-2/(d \vee 4)} \big),
\]
uniformly over $0 < \underline{\delta} \le 1$, where $d = d_x \wedge d_y$. The upper bound, which is attained by the empirical estimator, directly follows from the result of \cite{zhang2024gromov}. The contribution here is to supply a matching lower bound, up to logarithmic factors. 

We highlight key proof techniques. 
For $\Wp$, our lower-bound proofs use the method of fuzzy hypotheses (cf. Theorem 2.14 in \cite{tsybakov2009introduction}) and construct two pairs of \textit{random} (discrete) measures $(\mu_0,\nu_0), (\mu_1,\nu_1)$, which serve as priors, such that (a) they are supported on $\{ \Wp \ge \underline{\delta} \}$ and suitably separate $\Wp$, and that (b) the total variation between the marginal distributions of the observations is sufficiently small.\footnote{Condition (b) can be weakened without using the total variation; see Proposition \ref{prop: fuzzy} below.} To verify Condition (b), we use a reduction to testing for distributions over finite sets, similar to \cite{nilesweed2022estimation}. The bulk of our effort is devoted to constructing random measures satisfying Condition (a), which requires several new ideas. The construction differs between $1 \le p \le 2$ and $p > 2$. For $1\le p\le2$, we add isolated atoms to discrete measures and establish a decomposition of the transport cost that creates a prescribed separation (cf. Lemmas \ref{lem: Wp block} and \ref{lem: Wp exact}). For $p>2$, we instead use symmetric discrete measures and translate one marginal. Symmetry cancels the linear term in an expansion of the ground cost ($\|x-y\|^p$), and the remainder term is governed by the \textit{Bregman divergence}. Handling the remainder term is especially challenging for $4 < d \le 2p$ and entails carefully upper-bounding the transport cost for the random discrete measures for a \textit{mixed-power} ground cost.  

The upper bounds for the two regimes, $d >2p$ and $1 \le p \le 2$, and $d \le 4 \wedge (2p)$, follow directly from the existing results \cite{fournier2015on,manole2024sharp}. The proofs for the other regimes require several new ideas. Our proofs draw two independent samples from an optimal population coupling and match the midpoints of the sampled pairs. Averaging the costs of the two associated cross-couplings cancels the linear term in an expansion of the ground cost ($\|x-y\|^p$). For $p>2$, the remainder term is controlled by the Bregman divergence, leading to the factor $t^{p-2}$ in the empirical transport-cost bound (\ref{eq: transport-cost bound}) in the regime $d > 2p$. In the regime $4 < d \le 2p$, we use two ingredients in OT theory to derive sharp (in $\underline{\delta})$ bounds for the empirical estimator; one is a dyadic bound for the transport cost for a general continuous ground cost (Lemma \ref{lem: dyadic}), and the other is the \textit{shortening} principle for OT (Lemma \ref{lem: shortening}; cf.~Chapter 8 in \cite{villani2009optimal}). 

Turning to $\GW$, 
the lower-bound proof requires different prior constructions than $\Wp$ because of invariance of GW under isometries. Instead of adding atoms to random discrete measures, we dilate one marginal and use an identity that separates $\GW^2$ from $\underline{\delta}^2$ (cf. Lemma \ref{lem: GW-scaling}). We then use a concentration inequality from \cite{sambale2022concentration} to control a fluctuation of the remainder, which is random.

Altogether, the present paper contributes to a better understanding of the fundamental difficulty of estimating various OT distances under the vanishing separation regime, with several novel proof techniques. Our results provide sharp quantification of the benefit of distance separation for estimation of $\Wp$ and $\GW$, compared with unrestricted minimax rates.

\subsection{Literature review}

Establishing sharp rates of convergence for empirical distributions under $\Wp$ has been extensively studied in the statistics and probability literature; see, e.g., \cite{dudley1969the,ajtai1984on,talagrand1992matching,talagrand1994matching,talagrand1994the,dobric1995asymptotics,delbarrio1999central,boutetdemonvel2002almost,barthe2013combinatorial,dereich2013constructive,boissard2014on,fournier2015on,weed2019sharp,lei2020convergence,bobkov2019one,fournier2023convergence}. We shall recall Theorem 1 from \cite{fournier2015on}, which shows that, when $\mu$ is compactly supported, 
\[
\E[\Wp^p(\hat{\mu}_n,\mu)] \lesssim 
\begin{cases}
    n^{-1/2}, & d < 2p, \\
    n^{-1/2}\log n, & d = 2p, \\
    n^{-p/d}, & d > 2p.
\end{cases}
\]
Jensen's inequality yields $\E[\Wp(\hat{\mu}_n,\mu)] \le (\E[\Wp^p(\hat{\mu}_n,\mu)])^{1/p}$, and the triangle inequality yields (\ref{eq: baseline}) for $d > 2p$. 

Estimation of the transport cost $\Wp^p$, rather than distributional estimation under $\Wp$, has been explored in recent works by \cite{nilesweed2022estimation,chizat2020faster,manole2024sharp,hundrieser2024empirical,staudt2025convergence}. 
For $\Wp^p$ (and not $\Wp$), strong duality enables us to use empirical process techniques to bound $\E[|\Wp^p(\hat{\mu}_n,\hat{\nu}_n)-\Wp^p(\mu,\nu)|]$.
With the trivial inequality $|a-b| \le b^{1-p}|a^p-b^p|$ for $a \ge 0, b > 0, p > 1$, the results of \cite{chizat2020faster,manole2024sharp,hundrieser2024empirical} imply the rate (\ref{eq: known rate}) for $\Wp$. For $p = 2$ and compactly supported marginals, \cite{chizat2020faster} exploited semiconvexity and uniform Lipschitz estimates of dual potentials to establish $O(n^{-2/d})$ rates when $d > 4$. Extensions to other $p$ and unbounded marginals were taken up by \cite{manole2024sharp}; see also \cite{hundrieser2024empirical}. Those results provide sharp upper bounds for the empirical estimator with $\underline{\delta} \downarrow 0$ for the two cases where $d > 2p$ and $1 \le p \le 2$, and $d \le 4 \wedge (2p)$; see Table \ref{tab:separated-rates}. Our Theorems \ref{thm: Wp main} and \ref{thm: Wp main small p} complement the other cases by proving novel bounds for the empirical estimator under small $\underline{\delta}$. 

\begin{table}[htbp]
\centering
\caption{Minimax lower bounds and sharp empirical upper bounds under vanishing separation: Novelty comparison.}
\label{tab:separated-rates}
\small
\setlength{\tabcolsep}{3pt}
\renewcommand{\arraystretch}{1.13}
\begin{tabular*}{\linewidth}{@{\extracolsep{\fill}}llll@{}}
\toprule
 & & \multicolumn{2}{c@{}}{\bfseries Empirical upper bound} \\
\cmidrule(l){3-4}
\textbf{Setting} & \textbf{Minimax lower bound} & Case & Source \\
\midrule
$\Wp$, $d>2p$ & Theorem~\ref{thm: Wp main} & $1\le p\le2$
 & \cite{fournier2015on,hundrieser2024empirical} \\
 & & & \cite{manole2024sharp} ($d\ge5$), \cite{chizat2020faster} ($p=2$) \\
\addlinespace[2pt]
 & & $p>2$ & Theorem~\ref{thm: Wp main} \\
\midrule
$\Wp$, $d\le2p$ & Theorem~\ref{thm: Wp main small p} & $d\le4$
 & \cite{fournier2015on,hundrieser2024empirical} \\
 & & & \cite{chizat2020faster} ($p=2$) \\
\addlinespace[2pt]
 & & $d>4$ & Theorem~\ref{thm: Wp main small p} \\
\midrule
$\GW$, all $d$ & Theorem~\ref{thm: GW main} & & \cite{zhang2024gromov,kato2025convergence} \\
\bottomrule
\end{tabular*}
\end{table}

Regarding minimax lower bounds, the literature has focused on \textit{unrestricted} classes of distributions. For $\Wp$, Theorem 11 in \cite{nilesweed2022estimation} establishes, for $d > 2p$,
\[
\inf_{\hat{W}} \sup_{\mu,\nu} \E[|\hat{W}-\Wp(\mu,\nu)|] \gtrsim (n\log n)^{-1/d},
\]
where $\inf_{\hat{W}}$ is taken over all estimators based on $n$ independent observations from each of $\mu,\nu$, and $\sup_{\mu,\nu}$ is taken over all $\mu,\nu$ supported in a fixed ball. This result cannot translate into our lower-bound problem, simply because $\underline{\delta}$-separation restricts the class of distributions to $\{ (\mu,\nu) : \Wp(\mu,\nu) \ge \underline{\delta} \}$, and as such, the minimax risk under separation is smaller than or equal to that under the unrestricted class. Another related result is Theorem 22 in \cite{manole2024sharp}, which implies that, for $d > 2p$,
\begin{equation}
\inf_{\hat{T}} \sup_{\mu,\nu} \E[|\hat{T}-\Wp^p(\mu,\nu)|] \gtrsim (n\log n)^{-(p \wedge 2)/d},
\label{eq: minimax cost}
\end{equation}
where, again, $\sup_{\mu,\nu}$ is taken over all $\mu,\nu$ supported in a fixed ball.
We stress that Theorem 22 in \cite{manole2024sharp} does not yield sharp lower bounds 
in our problem when $\underline{\delta} \downarrow 0$, or for $1 < p < 2$ even when $\underline{\delta} > 0$ stays fixed; see Remark \ref{rem: manole} for details. The same remark partly explains a subtlety of our problem. The minimax lower bounds presented in Theorems \ref{thm: Wp main} and \ref{thm: Wp main small p} are thus new.  

The GW distance defines a metric on a space of Polish metric measure spaces modulo measure-preserving isometries \cite{memoli2011gromov,sturm2023the} and enables comparing and aligning heterogeneous and structured data sets. It has seen various applications, encompassing language alignment \cite{alvarezmelis2018gromov}, shape and graph matching \cite{memoli2009spectral, xu2019gromov, xu2019scalable}, single-cell genomics \cite{demetci2022scot}, to mention a few. 

For quadratic GW, the work \cite{zhang2024gromov} develops a useful variational representation of $\GW^2$, which enables establishing empirical convergence rates for compactly supported marginals. The work \cite{kato2025convergence} extends the sample complexity result to marginals with unbounded supports and establishes minimax lower bounds for $\GW^2$ (without separation), akin to (\ref{eq: minimax cost}), which, however, do not imply the lower bound in Theorem \ref{thm: GW main}. The very recent work \cite{leung2026rate}, posted on arXiv in November 2026,  obtains empirical convergence rates for general $(p,q)$-GW costs. Extending minimax lower bounds under separation for general $(p,q)$-GW distances will not be pursued in the present paper.

\subsection{Organization}
Section \ref{sec: preliminaries} introduces the two distances and the corresponding minimax risks. Section \ref{sec: main} states the main results. Sections \ref{sec: Wp proof}, \ref{sec: proof Wp small p}, and \ref{sec: GW proof} contain the proofs for Theorems \ref{thm: Wp main}, \ref{thm: Wp main small p}, and \ref{thm: GW main}, respectively. 

\subsection{Notation}

For any Polish metric space $\sfX$ endowed with its Borel $\sigma$-field, let $\calP(\sfX)$ denote the set of Borel probability measures. For $\mu \in \calP(\sfX)$ and a measurable mapping $F$ from $\sfX$ into another measurable space, we use $F_{\#}\mu$ to denote the pushforward, i.e., $F_{\#}\mu = \mu \circ F^{-1}$. For $a,b \in \R$, we use the notation $a \wedge b = \min \{a,b \}$ and $a \vee b = \max \{ a,b \}$. For any $k \in \N$, we write $[k] = \{ 1,\dots,k \}$. 
Let $B_d(x,r)$ denote the closed ball in $\R^d$ with center $x$ and radius $r>0$. We use $\lesssim_\xi$ to denote an inequality that holds up to a constant depending only on $\xi$; the subscript is dropped when the constant is universal. We write $a \asymp_\xi b$ if $a \lesssim_\xi b$ and $b \lesssim_\xi a$.

\section{Preliminaries}
\label{sec: preliminaries}
Let $\sfX$ and $\sfY$ be Polish spaces, and let $(\mu,\nu) \in \calP(\sfX) \times \calP(\sfY)$. 
Let $X_1,\dots,X_n$ be i.i.d. from $\mu$ and $Y_1,\dots,Y_m$ be i.i.d. from $\nu$, where the two samples are independent. We allow for unequal sample sizes, $n \ne m$. 
Let $N=n\wedge m$. The empirical distributions are
\[
 \hat\mu_n=\frac{1}{n}\sum_{i=1}^n\delta_{X_i}
 \quad\text{and}\quad
 \hat\nu_m=\frac{1}{m}\sum_{j=1}^m\delta_{Y_j}.
\]
In what follows, we only consider cases where $\sfX$ and $\sfY$ are finite-dimensional Euclidean spaces. 
\subsection{Optimal transport}
We first review OT. See, e.g., \cite{villani2009optimal,santambrogio2015optimal}, for details. 
For a continuous cost function $\sfc: \sfX \times \sfY \to \R_+$, the OT cost for $(\mu,\nu)$ is defined by 
\[
 T_{\sfc}(\mu,\nu):=\inf_{\pi \in \Pi(\mu,\nu)} \int \sfc \, d\pi, 
\]
where $\Pi(\mu,\nu)$ is the collection of couplings for $(\mu,\nu)$. Any coupling is a joint distribution with marginals $\mu,\nu$. 
The infimum on the right-hand side is attained, and the following strong duality holds
\begin{equation}
 T_\sfc(\mu,\nu)=\sup_{f,g}
\left ( \int f\,d\mu+\int g\,d\nu\right ),
 \label{eq: OT-dual}
\end{equation}
where the supremum is taken over all bounded continuous functions $f: \sfX\to \R$ and $g: \sfY \to \R$ such that $f(x) +g(y) \le \sfc(x,y)$ for all $(x,y)\in \sfX \times \sfY$. See Chapters 4 and 5 in \cite{villani2009optimal}. Strong duality yields joint convexity of the OT cost $T_{\sfc}$ in both marginals.

When $\sfX=\sfY=\R^d$, for $1 \le p<\infty$, the $p$-Wasserstein distance is defined by 
\[
\Wp(\mu,\nu):=T_\sfc(\mu,\nu)^{1/p} \quad \text{with} \quad \sfc(x,y) = \|x-y\|^p. 
\]
The $\Wp$ defines a distance on the subset of $\calP(\R^d)$ with finite $p$-th moments. In particular, $\Wp(\mu,\nu) = 0$ if and only if $\mu=\nu$. Convergence in $\Wp$ is equivalent to weak convergence plus convergence of $p$-th moments. 
When $(\mu,\nu) = (\hat{\mu}_n,\hat{\nu}_m)$, $\Wp^p(\hat\mu_n,\hat\nu_m)$ can be expressed as 
\begin{equation}
 \Wp^p(\hat\mu_n,\hat\nu_m)
 =\inf_{\gamma}\sum_{i=1}^n\sum_{j=1}^m\gamma_{ij}\|X_i-Y_j\|^p,
 \label{eq: discrete Wp}
\end{equation}
where the infimum is taken over all nonnegative $n \times m$ matrices $\gamma = (\gamma_{ij})$ satisfying
\[
\sum_{j=1}^m \gamma_{ij} = \frac{1}{n}, \quad \sum_{i=1}^n\gamma_{ij}=\frac{1}{m}. 
\]
When $n=m$, an optimal $\gamma$ is given by $1/n$ times a permutation matrix by the Birkhoff--von Neumann theorem \cite[Proposition~2.1]{peyre2019computational}. In particular, when $n=m$,
\begin{equation}
\Wp^p(\hat\mu_n,\hat\nu_n)
 =\inf_{\sigma}\frac{1}{n}\sum_{i=1}^n\|X_i-Y_{\sigma(i)}\|^p,
 \label{eq: matching}
\end{equation}
where the infimum is taken over all permutations $\sigma$ of $[n]$. The infimum on the right-hand side is attained. 

\subsection{Gromov--Wasserstein distance}
\label{sec: prelim-GW}
For Polish metric measure spaces $(\sfX,\mathsf{d}_{\sfX}, \mu),  (\sfY,\mathsf{d}_{\sfY},\nu)$, and for $p,q \in [1,\infty)$, the $(p,q)$-GW distance  is defined by
\begin{equation*}
\begin{split}
\mathsf{GW}_{p,q}^p(\mu,\nu) &:= \mathsf{GW}_{p,q}^p\Big( (\sfX,\mathsf{d}_{\sfX}, \mu), (\sfY,\mathsf{d}_{\sfY},\nu) \Big)  \\
&:= \inf_{\pi \in \Pi(\mu,\nu)} \iint \big | \mathsf{d}_{\sfX}(x,x')^q - \mathsf{d}_{\sfY}(y,y')^q \big |^p \, d\pi(x,y) d\pi(x',y').
\end{split}
\label{eq: GW}
\end{equation*}
 If one views $\mathsf{GW}_{p,q}$ as comparing two metric measure spaces,  then it is symmetric and satisfies the triangle inequality. Finally, $\mathsf{GW}_{p,q}(\mu,\nu) = 0$ if and only if there exists an isometry $f: \mathrm{spt} (\mu) \to \mathrm{spt}(\nu)$ such that $\nu = f_{\#}\mu$ ($\mathrm{spt}(\mu)$ denotes the support of $\mu$).  See \cite{memoli2011gromov,sturm2023the} for details. 
 
As in \cite{zhang2024gromov,kato2025convergence}, we shall focus on $\sfX = \R^{d_x}, \sfY = \R^{d_y}$, and $(p,q)=(2,2)$, where we endow $\R^{d_x}$ and $\R^{d_y}$ with the standard Euclidean metrics. Let $\mu \in \calP(\R^{d_x})$ and $\nu \in \calP(\R^{d_y})$ be compactly supported. Expanding the square, one can express the (squared) quadratic GW distance $\GW^2(\mu,\nu)$ as
\begin{equation}
\GW^2(\mu,\nu)=\sfA(\mu)+\sfA(\nu)
 -2\sup_{\pi\in\Pi(\mu,\nu)}
   \iint\|x-x'\|^2\|y-y'\|^2\,d\pi(x,y)\,d\pi(x',y'),
   \label{eq: GW-fourth-moments}
\end{equation}
where $\sfA(\mu) := \iint \|x-x'\|^4 \, d\mu(x)d\mu(x')$.

\subsection{Minimax risks}
For $\Wp$, we consider marginals supported in the (common) unit ball in $\R^d$. For $0 < \underline{\delta} \le 2$, define
\begin{equation*}
\label{eq: Wp risk}
 \calR^{\Wp}_{n,m,d}(\underline{\delta})
 :=\inf_{\hat{W}}\
   \sup_{\substack{\mu,\nu \in \calP(B_{d}(0,1))\\ \Wp(\mu,\nu)\ge\underline{\delta}}}
   \E_{\mu,\nu}\big[\bigl|\hat{W}-\Wp(\mu,\nu)\bigr|\big].
\end{equation*}

For $\GW$, we consider marginals supported in the unit balls of their respective spaces. For $0 <  \underline{\delta} \le 2\sqrt{2}$, define
\begin{equation*}
\label{eq: GWp risk}
 \calR^{\mathsf{GW}}_{n,m,d_x,d_y}(\underline{\delta})
 :=\inf_{\hat G}\
   \sup_{\substack{\mu\in\calP(B_{d_x}(0,1)),\nu\in\calP(B_{d_y}(0,1))\\
                    \GW(\mu,\nu)\ge\underline{\delta}}}
   \E_{\mu,\nu}\big[\bigl|\hat G-\GW(\mu,\nu)\bigr|\big].
\end{equation*}
In both cases, the infimum is taken over measurable functions of the two samples, and $\E_{\mu,\nu}[\cdot]$ denotes expectation under $\mu^{\otimes n}\otimes\nu^{\otimes m}$. The empirical estimators are $\Wp(\hat\mu_n,\hat\nu_m)$ and $\GW(\hat \mu_n,\hat \nu_m)$, respectively.

\begin{remark}[Nonemptiness of separation classes]
\label{rem: nonemptiness}
The constraints $\underline{\delta} \le 2$ and $\underline{\delta} \le 2\sqrt{2}$ guarantee that the classes of marginals $(\mu,\nu)$ for which $\Wp(\mu,\nu) \ge \underline{\delta}$ and $\GW(\mu,\nu) \ge \underline{\delta}$ are nonempty. To see this, for a unit vector $v \in \R^d$ with $d_x=d$ and $0 \le t \le 1$, we have $\Wp(\delta_{tv},\delta_{-tv}) = 2t$ and $\GW (\mu_t,\delta_0) = 2\sqrt{2}t^2$ for $\mu_t=(\delta_{tv}+\delta_{-tv})/2$. These are rather trivial distributions. If one restricts to absolutely continuous distributions on the unit balls, both for $(\mu,\nu)$, then the range of $\Wp$ is exactly $[0,2)$, and the range of $\GW$ is exactly $[0,2\sqrt{2})$ if $d_x = d_y$ and $(0,2\sqrt{2})$ if $d_x \ne d_y$; see Lemma \ref{lem: nonemptiness}. 
Conversely, it is obvious that these classes are empty if $\underline{\delta} > 2$ for $\Wp$ and if $\underline{\delta} > 2\sqrt{2}$ for $\GW$. 
\end{remark}

We collect a few immediate properties of the minimax risks. 

\begin{lemma}
\label{lem: minimax risk}
Consider the minimax risks $\calR^{\mathsf{\Wp}}_{n,m,d}(\underline{\delta})$ and $\calR^{\mathsf{GW}}_{n,m,d_x,d_y}(\underline{\delta})$.
\begin{enumerate}
\item[(i)] Both $\calR^{\mathsf{\Wp}}_{n,m,d}(\underline{\delta})$ and $\calR^{\mathsf{GW}}_{n,m,d_x,d_y}(\underline{\delta})$ are nonincreasing in $\underline{\delta}$. 
\item[(ii)] For $d = d_x \wedge d_y$, $\calR^{\mathsf{GW}}_{n,m,d_x,d_y}(\underline{\delta}) \ge \calR^{\mathsf{GW}}_{n,m,d,d}(\underline{\delta})$.
\item[(iii)] For $d > d'$, $\calR^{\mathsf{\Wp}}_{n,m,d}(\underline{\delta}) \ge \calR^{\mathsf{\Wp}}_{n,m,d'}(\underline{\delta})$ and $\calR^{\mathsf{GW}}_{n,m,d,d}(\underline{\delta}) \ge \calR^{\mathsf{GW}}_{n,m,d',d'}(\underline{\delta})$.
\end{enumerate}
\end{lemma}
\begin{proof}
Part (i) is trivial. For Part (ii), if $d_y \ge d_x$, then
by considering $\nu$ such that the last $d_y-d_x$ coordinates of $Y \sim \nu$ are degenerate to $0$, we see that $\calR^{\mathsf{GW}}_{n,m,d_x,d_y} \ge \calR^{\mathsf{GW}}_{n,m,d,d}$. The same inequality holds when $d_x \ge d_y$. Part (iii) follows from the fact that the pushforward by $\R^{d'} \ \ni x \mapsto (x,0) \in \R^{d}$ is isometric from $\calP(B_{d'}(0,1))$ into $\calP(B_d(0,1))$, with respect to $\Wp$ and $\GW$.
\end{proof}

\section{Main results}
\label{sec: main}
\subsection{Wasserstein distance}
\label{sec: Wp}

We first establish minimax rates for $\Wp$ when $\Wp \ge \underline{\delta}$ for some $\underline{\delta} > 0$. The next theorem concerns the $d > 2p$ case. 
\begin{theorem}[Minimax rates for $\Wp$: $d > 2p$ case] 
\label{thm: Wp main}
Let $p \ge 1, d>2p$, and $\alpha=p\wedge2$.  Put  $\ell_N=(N\log N)^{-1/d}$ for $N \ge 2$. Then, there exists a positive constant $N_0$, depending only on $p,d$, such that for $N\ge N_0$ and $0<\underline{\delta} \le 3/4$,
\begin{equation*}
\label{eq: Wp main}
\ell_N\wedge(\underline{\delta}^{1-\alpha}\ell_N^\alpha) \lesssim_{p,d}
  \calR^{\Wp}_{n,m,d}(\underline{\delta})
 \lesssim_{p,d} N^{-1/d}\wedge(\underline{\delta}^{1-\alpha}N^{-\alpha/d}).
\end{equation*}
The empirical estimator $\Wp(\hat\mu_n,\hat\nu_m)$ attains the upper bound. 
\end{theorem}

The theorem shows that the minimax rate agrees with 
\[
N^{-1/d} \wedge (\underline{\delta}^{1-\alpha} N^{-\alpha/d})
\]
up to logarithmic factors, where $\underline{\delta} = \underline{\delta}_N \downarrow 0$ is allowed. Assume $p > 1$. When $\underline{\delta}=N^{-\beta/d}$ for some $\beta > 0$, 
\[
N^{-1/d} \wedge (\underline{\delta}^{1-\alpha} N^{-\alpha/d}) = N^{-1/d} \wedge N^{-(\alpha - (\alpha-1)\beta)/d} = 
\begin{cases}
N^{-1/d} & \text{if $\beta > 1$}, \\
N^{-(\alpha - (\alpha-1)\beta)/d} & \text{if $\beta \le 1$}.
\end{cases}
\]
In particular, when $p \ge 2$ (i.e., $\alpha = 2$), the rate reduces to $N^{-1/d}$ for $\beta > 1$ and $N^{-(2-\beta)/d}$ for $\beta \le 1$, which interpolates between $N^{-1/d}$ and $N^{-2/d}$.

For $1 \le p \le 2$, the upper bound follows directly from the results of \cite{chizat2020faster, manole2024sharp, hundrieser2024empirical}. For $p > 2$, Theorem \ref{thm: Wp main} improves the dependence of the empirical upper bound on $\underline{\delta}$, from $\underline{\delta}^{1-p}$ to $\underline{\delta}^{-1}$, which follows by a novel distance-dependent bound for the empirical transport cost; see Proposition \ref{prop: power bound} below. The improvement can be polynomial in $N$ when $\underline{\delta}$ is vanishing. For example, with $p=3, d=8$, and $\underline{\delta}=N^{-1/16}$,  the preceding bound (\ref{eq: known rate}) gives $N^{-1/8}$, whereas our bound gives $N^{-3/16}$.

\begin{remark}[Comparison with Theorem 22 in \cite{manole2024sharp}]
\label{rem: manole}
For $d >2p >2$, Theorem 22 in \cite{manole2024sharp} implies a lower bound of order $(N \log N)^{-2/d}$ for $\calR^{\Wp}_{n,m,d}(\underline{\delta})$ for sufficiently small $\underline{\delta}$. The exponent here is $2/d$ and not $\alpha/d$, to be clarified below. Assume $n=m$.  In their notation, take (sufficiently small but fixed) balls $\calX, \calY \subset B_d(0,1)$ that are separated by at least $c > 0$, independent of $n$. Observe that Theorem 22 in \cite{manole2024sharp} does not allow $\calX, \calY$ to depend on $n$. For $(\mu,\nu) \in \calP(\calX) \times \calP(\calY)$, $\Wp(\mu,\nu) \ge c$. Since $\Wp \le 2$ on $\calP(\calX) \times \calP(\calY)$, it suffices to consider estimators $\hat{W}$ with values in $[0,2]$. For such estimators, we have
\[
|\hat{W}^p - \Wp^p| \le p2^{p-1}|\hat{W} -\Wp|.
\]
For separated supports, the assumptions of Theorem 22 in \cite{manole2024sharp} hold with $\alpha=2$, including when $1 < p < 2$ (observe that the $p$-th power function is smooth on any open set not containing the origin). 
Since $\{ (\mu,\nu) \in \calP(B_d(0,1)) \times \calP(B_d(0,1)) : \Wp(\mu,\nu) \ge \underline{\delta} \} \supset  \calP(\calX) \times \calP(\calY)$, Theorem 22 in \cite{manole2024sharp} implies, for $d > 2p >2$ and $0 < \underline{\delta} \le c$,
\[
\begin{split}
\calR^{\Wp}_{n,n,d} (\underline{\delta}) &\ge \inf_{\hat{W} \in [0,2]} \sup_{(\mu,\nu) \in \calP(\calX) \times \calP(\calY)} \E_{\mu,\nu}[|\hat{W} - \Wp|] \\
&\gtrsim_p \inf_{\hat{T}} \sup_{(\mu,\nu) \in \calP(\calX) \times \calP(\calY)} \E_{\mu,\nu}[|\hat{T} - \Wp^p|] \\
&\gtrsim_{p,d} (n \log n)^{-2/d},
\end{split}
\]
which is not tight  when $\underline{\delta} = \underline{\delta}_n \downarrow 0$.  Even when $\underline{\delta} > 0$ stays fixed, the order $(n \log n)^{-2/d}$ is not tight for $1 < p < 2$. One might predict that, for $p > 1$, by taking \textit{shrinking} balls $\calX, \calY$ with radii $\underline{\delta}$ that are separated by $\underline{\delta}$, a lower bound of order $\underline{\delta}^{1-p}(n \log n)^{-2/d}$ would follow, which, however,  would contradict the upper bound if $\underline{\delta} = \underline{\delta}_n \downarrow 0$ fast enough that $n^{-1/d} = o(\underline{\delta}^{1-p}(n \log n)^{-2/d})$. To the best of the authors' knowledge, a simple reduction to Theorem 22 in \cite{manole2024sharp} would not lead to the lower bound in Theorem \ref{thm: Wp main} for $d > 2p$. 
\end{remark}

Next, we consider the $d \le 2p$ case.

\begin{theorem}[Minimax rates for $\Wp$: $d \le 2p$ case]
\label{thm: Wp main small p}
Let $p \ge 1$ and $1 \le d \le 2p$.
Put 
\begin{equation*}
\xi_N (\underline{\delta})= 
 N^{-1/(2p)}\wedge
 \begin{cases}
 \underline{\delta}^{1-p}N^{-1/2}, & d\le4,\\
 \big( \underline{\delta}^{1-p}N^{-1/2} +\underline{\delta}^{-1}N^{-2/d}\big),& 4<d,\ p\le d,\\
 \big( \underline{\delta}^{1-p}N^{-1/2}+\underline{\delta}^{1-2p/d}N^{-2/d} \big), & 4<d<p.
 \end{cases}
 \label{eq: Wp rate small p}
\end{equation*}
Then, there exists a positive constant $N_0$, depending only on $p,d$,
such that for $N\ge N_0$ and $0<\underline{\delta}\le 3/4$, 
\[
(\log N)^{-1} \xi_{N}(\underline{\delta}) \lesssim_{p,d} \calR^{\Wp}_{n,m,d} (\underline{\delta}) \lesssim_{p,d} (\log N) \xi_N(\underline{\delta}). 
\]
The empirical estimator $\Wp(\hat\mu_n,\hat\nu_m)$ attains the upper
bound. The logarithmic factor may be omitted from the lower bound for
$d\le4$ and from both bounds for $d<\min\{4,2p\}$.
\end{theorem}

For $d \le 4$, the lower bound follows by a simple application of Le Cam's two-point method. As discussed in the introduction, the lower- and upper-bound proofs for the complementary case require substantive effort; see the discussion at the start of Section \ref{sec: Wp small p prelim}. 

For $d > 4$, there is a phase transition in the dependence of $\underline{\delta}$ at $p=d$.  For $4 < d < p$, our lower-bound proof constructs random discrete measures and takes mixtures with point masses with mixture weight $\underline{\delta}^p$, which guarantees a prescribed separation. The ``effective sample'' size for the random discrete measure component now reduces to $O(N\underline{\delta}^p)$ (see (\ref{eq: sample size choice}) below), which translates into $\delta^{1-2p/d}N^{-2/d}$ in the lower bound.

\subsection{Quadratic Gromov--Wasserstein distance}
\label{sec: GW}
Next, we consider $\GW$.

\begin{theorem}[Minimax rates for $\GW$]
\label{thm: GW main}
For $N = n \wedge m \ge 2$ and $d = d_x \wedge d_y$, let
\begin{equation*}
\underline{\varphi}_N = N^{-1/4} \vee (N\log N)^{-1/d} \quad \text{and} \quad \overline{\varphi}_N = N^{-1/(d \vee 4)} \sqrt{(\log N)^{\ind_{\{ d=4 \}}}}. 
 \label{eq: GW-rho}
\end{equation*}
There exists a positive constant $N_0$, depending only on $d_x,d_y$, such that, for $N\ge N_0$ and $0<\underline{\delta} \le 1$,
\begin{equation}\label{eq: GW-main}
\underline{\varphi}_N \wedge (\underline{\delta}^{-1} \underline{\varphi}_N^2) \lesssim_{d_x,d_y}
   \calR^{\mathsf{GW}}_{n,m,d_x,d_y}(\underline{\delta})
  \lesssim_{d_x,d_y} \overline{\varphi}_N \wedge (\underline{\delta}^{-1} \overline{\varphi}_N^2).
\end{equation}
The empirical estimator $\GW(\hat\mu_n,\hat\nu_m)$ attains the upper bound. 
\end{theorem}

The upper bound follows directly by \cite{zhang2024gromov}. The lower bound is new. Up to logarithmic factors, the minimax rate agrees with
\[
N^{-1/(d \vee 4)} \wedge (\underline{\delta}^{-1}N^{-2/(d \vee 4)}),
\]
where the smaller of two dimensions, $d = d_x \wedge d_y$, determines the rate. A minimax lower bound for $\GW^2$ without separation was established by \cite{kato2025convergence}. For $\GW$, because of invariance under isometries, constructing random measures satisfying a prescribed separation requires new ideas, compared with the $\Wp$ case as well as \cite{kato2025convergence}. To this end, we use dilating marginals to create a separation; see Lemma \ref{lem: GW-scaling}.

\section{Proofs of Theorem \ref{thm: Wp main}}
\label{sec: Wp proof}

\subsection{Method of two fuzzy hypotheses}
Our lower-bound proofs rest on the method of two fuzzy hypotheses (cf. Theorem 2.14 in \cite{tsybakov2009introduction}), adapted to our use case.

Let $\Theta$ be a measurable space, and let $T: \Theta \to \R$ be a measurable function. For a Polish  space $\sfX$ endowed with the Borel $\sigma$-field, consider a statistical model $\{ P_\theta : \theta \in \Theta \}$, where $\theta \mapsto P_\theta$ is measurable from $\Theta$ into $\calP(\sfX)$. The latter space is endowed with the Borel $\sigma$-field corresponding to the topology of weak convergence. For any (measurable) estimator $\hat{T}: \sfX \to \R$, we are interested in finding a lower bound on the maximum risk
\[
\sup_{\theta \in \Theta}\E_\theta [|\hat{T}-T(\theta)|]. 
\]
Markov's inequality implies that, for $s > 0$, 
\[
\sup_{\theta \in \Theta}\E_\theta [|\hat{T}-T(\theta)|] \ge s \sup_{\theta \in \Theta}P_\theta (|\hat{T}-T(\theta)| \ge s).
\]
Let $\Lambda_0,\Lambda_1$ be two ``prior'' distributions on $\Theta$, and let $\overline{P}_0, \overline{P}_1$ be the corresponding marginal distributions on $\sfX$,
\[
\overline{P}_j (\cdot) = \int_\Theta P_\theta (\cdot) \, d\Lambda_j(\theta), \ j=0,1.
\]

\begin{proposition}[A version of Theorem 2.14 in \cite{tsybakov2009introduction}]
\label{prop: fuzzy}
Suppose that:
\begin{enumerate}
\item[(i)] 
There exist $c \in \R, s > 0$, and $0 \le \beta_0,\beta_1 < 1$ such that 
\[
\begin{split}
&\Lambda_0 (\{ \theta : T(\theta) \le c \}) \ge 1-\beta_0, \\
&\Lambda_1 (\{ \theta : T(\theta) \ge c + 2s \}) \ge 1-\beta_1.
\end{split}
\]
\item[(ii)] There exists $0 < \alpha < 1$ such that
\[
\overline{P}_0(\hat{T} \ge c+s ) +\overline{P}_1(\hat{T} < c+s ) \ge 1-\alpha. 
\]
\end{enumerate}
Then, we have 
\[
\sup_{\theta \in \Theta} P_\theta (|\hat{T}-T(\theta)| \ge s) \ge \frac{1-\alpha - \beta_0-\beta_1}{2}. 
\]
\end{proposition}

\begin{proof}
This follows directly from the proof of Theorem 2.14 in \cite{tsybakov2009introduction}. 
\end{proof}

\subsection{Preliminaries}

We will often use the following estimate on the $p$-th root. 
\begin{lemma}
\label{lem: rootgap}
For $p\ge1$, $0\le a<b$, $\delta\ge0$, and $c>0$,
\begin{equation}
  (\delta^p+bc)^{1/p}-(\delta^p+ac)^{1/p}
  \asymp_{p, a, b} c^{1/p}\wedge(\delta^{1-p}c).
       \label{eq: rootgap}
\end{equation}
At $\delta=0$, take the right-hand side of (\ref{eq: rootgap}) to be $c^{1/p}$.
\end{lemma}

\begin{proof}
Put $v=(\delta^p+ac)^{1/p}$ and $w=(\delta^p+bc)^{1/p}$. For $u=v/w \in [0,1)$, we have $u^p \le u$, and convexity of $x\mapsto x^p$ at $1$ gives $u^p\ge1+p(u-1)$. Multiplying $1-u \le 1-u^p\le p(1-u)$ by $w^p$ yields
\[
  w^{p-1}(w-v) \le w^p-v^p=(b-a)c\le p w^{p-1}(w-v).
\]
Divide by $w^{p-1}$ and use $w^p=\delta^p+bc\asymp_b\delta^p+c\asymp\delta^p\vee c$ to obtain (\ref{eq: rootgap}).
\end{proof}

Similar to \cite{nilesweed2022estimation}, our lower-bound proofs use a reduction to a testing problem for distributions over a finite set. We shall recall two technical propositions from \cite{nilesweed2022estimation}. 

For probability measures $P,Q$ defined on a common measurable space, recall that the total variation $\TV(P,Q)$ and the $\chi^2$-divergence $\chi^2(P,Q)$ are defined by
\begin{itemize}
\item (Total variation) 
\[
\TV(P,Q) := \sup_{A} |P(A)-Q(A)|;
\]
\item ($\chi^2$-divergence) 
\[
\chi^2(P,Q):=
\begin{cases}
    \int \left ( \frac{dP}{dQ} - 1\right )^2 \, dQ & \text{if $P \ll Q$}, \\
    \infty & \text{otherwise}.
\end{cases}
\]
\end{itemize}
If $P,Q \ll R$ for some $\sigma$-finite measure $R$, then Scheff\'{e}'s lemma yields \begin{equation}
\TV(P,Q) = 1-\int \min \left \{ \frac{dP}{dR}, \frac{dQ}{dR} \right\} \, dR = \frac{1}{2} \int \left |  \frac{dP}{dR} - \frac{dQ}{dR} \right | \, dR.
\label{eq: scheffe}
\end{equation}
Applying the Cauchy--Schwarz inequality gives 
\[
\TV(P,Q) \le \frac{1}{2} \sqrt{\chi^2(P,Q)}.
\]
See Chapter 2 in \cite{tsybakov2009introduction}. 

\begin{lemma}[Proposition 9 in \cite{nilesweed2022estimation}]
\label{lem: Wp packing}
Let $\calB$ be a ball in $\R^d$ with radius $r > 0$. For any integer $M\ge2$, there exist $M$ points in $\mathcal{B}$
at pairwise distances at least $c_0 M^{-1/d}$, where $c_0$ depends
only on $d$ and $r$. Let $\mathcal{S}$ be the set of these points, and choose
$F$ uniformly among the bijections from $[M]$ onto $\mathcal{S}$.
Let $\mathfrak{u}$ be the uniform distribution on $[M]$. For any fixed $\mathfrak{q} \in \calP([M])$, provided that $d > 2p$, it holds that
\begin{equation}
 c_0^p M^{-p/d}\TV(\mathfrak{q},\mathfrak{u}) 
 \le \Wp^p(F_{\#}\mathfrak{q},F_{\#}\mathfrak{u}) 
 \le C_0 M^{-p/d}\bigl(\chi^2(\mathfrak{q},\mathfrak{u})\bigr)^{p/d}\TV(\mathfrak{q},\mathfrak{u})^{1-2p/d}
\label{eq: packing}
\end{equation}
with probability at least $0.9$, where $C_0 > 0$ depends only on $p, d$, and $r$.
\end{lemma}

 For $0<\eta\le1/10$, define
\begin{align*}
  \Dminus&:=\{\mathfrak{q} \in \calP([M]):\chi^2(\mathfrak{q},\mathfrak{u})\le9,\ \TV(\mathfrak{q},\mathfrak{u})\le\eta\},\\
  \Dplus&:=\{\mathfrak{q} \in \calP([M]):\chi^2(\mathfrak{q},\mathfrak{u})\le9,\ \TV(\mathfrak{q},\mathfrak{u})\ge1/4\}.
\end{align*}

\begin{lemma}[Proposition 10 in \cite{nilesweed2022estimation}]
\label{lem: Wp testing}
Fix $\eta\in(0,1/10]$ and let $k$ be the sample size. For a sufficiently large universal constant $C_1$, set
\begin{equation}
\label{eq: M-choice}
  M=\lceil C_1\eta^{-1}k\log k\rceil.
\end{equation}
For all sufficiently large $k$ (depending on $\eta$), any (possibly randomized) test $\psi$ with values in $\{ 0,1 \}$, based on $k$ independent observations from $\mathfrak{q}$, satisfies
\begin{equation*}
\label{eq: Wp testing-bound}
  \sup_{\mathfrak{q} \in \Dminus}\Prob_{\mathfrak{q}}(\psi=1)+
  \sup_{\mathfrak{q} \in \Dplus}\Prob_{\mathfrak{q}}(\psi=0)\ge 0.9.
\end{equation*}
\end{lemma}

\subsection{Lower bound for $1 \le p \le 2$}
\label{sec: lower bound small p}

By symmetry, we may assume $n \le m$, i.e., $N=n$. For ease of notation, set $\calX = B_d(0,1)$. We shall apply Proposition \ref{prop: fuzzy} with $\sfX=\calX^{n+m}, \Theta = \{ (\mu,\nu) \in \calP(\calX) \times \calP(\calX): \Wp(\mu,\nu) \ge \underline{\delta} \}, P_\theta = \mu^{\otimes n} \otimes \nu^{\otimes m}$, and $T(\theta) = \Wp(\mu,\nu)$ for $\theta = (\mu,\nu)$. To construct priors, we consider a random discrete measure of the form $F_{\#}\mathfrak{q}$ for some $\mathfrak{q} \in \calP([M])$ and the random bijection $F: [M] \to \calS$ from Lemma \ref{lem: Wp packing}, for a sufficiently small ball $\calB \subset \calX$, and add an atom to $F_{\#}\mathfrak{q}$  to construct $\mu_{F,\mathfrak{q}}$. Combined with a fixed marginal $\overline{\nu}$,  we shall construct such $\mu_{F,\mathfrak{q}}$ so that:
\begin{itemize}
    \item $(\mu_{F,\mathfrak{q}},\overline{\nu})$ is supported in $\{ \Wp \ge \underline{\delta} \}$; and
    \item $\Wp (\mu_{F,\mathfrak{q}^{-}},\overline{\nu})$ for $\mathfrak{q}^{-} \in \Dminus$ and $\Wp(\mu_{F,q^{+}},\overline{\nu})$ for $\mathfrak{q}^{+} \in \Dplus$ are separated by $\ell_N \wedge (\underline{\delta}^{1-p}\ell_N^p)$ up to a constant; cf. Condition (i) in Proposition \ref{prop: fuzzy}.
\end{itemize}
Once such random discrete measures are constructed, one can resort to Lemma \ref{lem: Wp testing} to verify Condition (ii) in Proposition \ref{prop: fuzzy}.
Constructing such $\mu_{F,\mathfrak{q}}$ requires delicate care. The next lemma will play a key role.  

\begin{lemma}\label{lem: Wp block}
Let $\mu_i\in\calP(\mathcal{X}_i)$ and $\nu_i\in\calP(\mathcal{Y}_i)$ for $i\in\{0,1\}$, where
the four sets are compact subsets of $\R^d$ and $\mathcal{X}_0\cap \mathcal{X}_1=\mathcal{Y}_0\cap \mathcal{Y}_1=\varnothing$.
Let $\overline{m}_{ij}$ and $\underline{m}_{ij}$ be the supremum and infimum of $\|x-y\|^p$
over $\mathcal{X}_i\times \mathcal{Y}_j$, respectively. For $0<\vartheta<1$ put
$\mu=\vartheta\mu_0+(1-\vartheta)\mu_1$ and $\nu=\vartheta\nu_0+(1-\vartheta)\nu_1$.
If $\overline{m}_{00}+\overline{m}_{11}<\underline{m}_{01}+\underline{m}_{10}$, then any optimal coupling for $\Wp(\mu,\nu)$ charges zero
mass to $(\mathcal{X}_0\times \mathcal{Y}_1)\cup(\mathcal{X}_1\times \mathcal{Y}_0)$, and
\begin{equation}
 \Wp^p(\mu,\nu)=\vartheta \Wp^p(\mu_0,\nu_0)+(1-\vartheta)\Wp^p(\mu_1,\nu_1).
\label{eq: coupling}
\end{equation}
\end{lemma}

\begin{proof}
Let $\pi$ be an optimal coupling for $\Wp(\mu,\nu)$ and put $\sfc(x,y)=\|x-y\|^p$. Since $\mu(\mathcal X_0)=\nu(\mathcal Y_0)=\vartheta$,
\[
 \pi(\mathcal X_0\times\mathcal Y_1)
 =\vartheta-\pi(\mathcal X_0\times\mathcal Y_0)
 =\pi(\mathcal X_1\times\mathcal Y_0).
\]
By Theorem 5.10~(ii) in \cite{villani2009optimal}, $\pi$ is concentrated on a $\sfc$-cyclically monotone set $\Gamma$. If the common mass were positive, $\Gamma$ would contain $(x_0,y_1)$ and $(x_1,y_0)$ with $x_i\in\mathcal X_i$ and $y_i\in\mathcal Y_i$ for $i\in\{0,1\}$. Cyclical monotonicity would imply
\[
 \underline{m}_{01}+\underline{m}_{10}\le \sfc(x_0,y_1)+\sfc(x_1,y_0)
 \le \sfc(x_0,y_0)+\sfc(x_1,y_1)\le \overline{m}_{00}+\overline{m}_{11},
\]
which contradicts the assumption. As such, $\pi(\mathcal X_0\times\mathcal Y_1) = \pi(\mathcal X_1\times\mathcal Y_0) = 0$, which guarantees
\[
\pi_0 := \vartheta^{-1} \pi|_{\calX_0 \times \calY_0} \in \Pi(\mu_0,\nu_0) \quad \text{and} \quad \pi_1 := (1-\vartheta)^{-1}\pi|_{\calX_1 \times \calY_1} \in \Pi(\mu_1,\nu_1). 
\]
Since $\pi = \vartheta \pi_0 + (1-\vartheta) \pi_1$, we have 
\begin{equation}
\Wp^p(\mu,\nu) = \int \sfc \, d\pi = \vartheta \int \sfc \, d\pi_0 + (1-\vartheta) \int \sfc \, d\pi_1 \ge \vartheta \Wp^p(\mu_0,\nu_0) + (1-\vartheta) \Wp^p(\mu_1,\nu_1). 
\label{eq: coupling 1}
\end{equation}

On the other hand, let $\pi_0',\pi_1'$ be optimal couplings for $\Wp(\mu_0,\nu_0), \Wp(\mu_1,\nu_1)$, respectively. Then $\pi' := \vartheta \pi_0' + (1-\vartheta)\pi_1'$ is a coupling for $(\mu,\nu)$, and 
\begin{equation}
\Wp^p(\mu,\nu) \le \int \sfc \, d\pi' = \vartheta \int \sfc \, d\pi_0' + (1-\vartheta) \int \sfc \, d\pi_1' = \vartheta \Wp^p(\mu_0,\nu_0) + (1-\vartheta) \Wp^p(\mu_1,\nu_1).
\label{eq: coupling 2}
\end{equation}
Combining (\ref{eq: coupling 1}) and (\ref{eq: coupling 2}) yields (\ref{eq: coupling}). 
\end{proof}

Let $e_1$ be the first coordinate vector. For a given positive integer $M \ge 2$ to be specified below, apply Lemma \ref{lem: Wp packing} with $\calB = B_d(7e_1/8,1/8)$ to construct an $M$-point set $\calS$ and the random bijection $F:[M]\to\mathcal S$. For $0<\underline{\delta}\le1$, set
\[
z_1=-e_1, \ z_2=e_1/2, \ \vartheta=(2\underline{\delta}/3)^p.
\]
Then $1-\vartheta\ge1/3$. Recall that $\mathfrak{u}$ is the uniform distribution on $[M]$. Pick any $\mathfrak{q} \in \calP([M])$. Now, consider
\[
 \mu_{F,\mathfrak{q}}=\vartheta\delta_{z_1}+(1-\vartheta)F_{\#}\mathfrak{q}
 \quad \text{and} \quad
 \overline{\nu}=\vartheta\delta_{z_2}+(1-\vartheta)F_{\#}\mathfrak{u} .
\]
Observe that, since $F_{\#}\mathfrak{u}$ is the uniform distribution on $\calS$, $\overline{\nu}$ does not depend on $F$. 
The map $F\mapsto\mu_{F,\mathfrak{q}}$ is measurable from $[M]^M$ into $\calP(\calX)$. We denote 
\begin{equation}
\Lambda_{\mathfrak{q}} = \text{the distribution of $(\mu_{F,\mathfrak{q}},\overline{\nu})$ on $\calP(\calX) \times \calP(\calX)$}. 
\label{eq: prior}
\end{equation}

Given $F$, we generate $(X_1,\dots,X_n,Y_1,\dots,Y_m)$ independently as $X_i \sim \mu_{F,\mathfrak{q}}$ and $Y_j \sim \overline{\nu}$. Concretely, one can  generate $(X_1,\dots,X_n,Y_1,\dots,Y_m)$ as follows. 
\begin{itemize}
    \item Generate $F,J_1,\dots,J_n, I_1,\dots,I_n$ independently such that $\Prob(J_i = j) = \vartheta^j (1-\vartheta)^{1-j}$ for $j \in \{ 0,1 \}$ and $I_i \sim \mathfrak{q}$. 
    \item Generate $X_1,\dots,X_n$ as $X_i = J_i z_1 + (1-J_i) F(I_i)$.
    \item Finally, generate $Y_1,\dots,Y_m \stackrel{\text{i.i.d.}}{\sim} \overline{\nu}$ independently from everything else. 
\end{itemize}

The next lemma follows from Lemma \ref{lem: Wp block} and our construction of $\mu_{F,\mathfrak{q}}$.

\begin{lemma}
\label{lem: Wp exact}
For any 
$\mathfrak{q} \in \calP([M])$,  $\mu_{F,\mathfrak{q}}$ and $\overline{\nu}$ are supported in $\calX$ and satisfy
\begin{equation}
 \Wp^p(\mu_{F,\mathfrak{q}},\overline{\nu})=\underline{\delta}^p+(1-\vartheta)\Wp^p(F_{\#}\mathfrak{q},F_{\#}\mathfrak{u}).
 \label{eq: Wp decomp}
\end{equation}
\end{lemma}

\begin{proof}
Take $\mathcal X_0=\{z_1\}$, $\mathcal Y_0=\{z_2\}$, and $\mathcal X_1=\mathcal Y_1=B_d(7e_1/8,1/8)$ in Lemma \ref{lem: Wp block}. The atoms $z_1$ and $z_2$ are $3/2$ apart, and the ball $B_d(7e_1/8,1/8)$ has diameter $1/4$. Since every point in $B_d(7e_1/8,1/8)$ has first coordinate at least $3/4$, its distances to $z_1$ and $z_2$ are at least $7/4$ and $1/4$, respectively. Hence
\[
\overline{m}_{00}+\overline{m}_{11} \le (3/2)^p+(1/4)^p < (7/4)^p+(1/4)^p \le \underline{m}_{01}+\underline{m}_{10}.
\]
Since $\vartheta\|z_1-z_2\|^p=\underline{\delta}^p$, we obtain (\ref{eq: Wp decomp}) from (\ref{eq: coupling}).
\end{proof}

We are now in position to prove the lower bound for $1 \le p \le 2$. 

\begin{proof}[Proof of Theorem \ref{thm: Wp main}: Lower bound for $1 \le p \le 2$]
Recall that we have assumed that $n \le m$, so that $N=n$. 
We follow the construction before Lemma \ref{lem: Wp exact}. 

Take $c_0$ and $C_0$ from Lemma \ref{lem: Wp packing} for $\calB = B_d(7e_1/8,1/8)$. 
Since $1-2p/d>0$, one may choose sufficiently small $0<\eta\le1/10$ (that may depend on $(p,d)$) such that
\[
C_0 9^{p/d}\eta^{1-2p/d} \le \frac{c_0^p}{16}.
\]
Choose $M$ as in (\ref{eq: M-choice}) with $k=n$ and set
\[
\Delta=\frac{(1-\vartheta)c_0^p M^{-p/d}}{16} \asymp_{p,d} (n \log n)^{-p/d}.
\]
Observe from (\ref{eq: Wp decomp}) that $\Wp^p(\mu_{F,\mathfrak{q}},\overline{\nu}) \ge \underline{\delta}^p$ for any $\mathfrak{q} \in \calP([M])$. 
Pick any $\mathfrak{q} \in \calP([M])$.

Suppose that the inequalities in (\ref{eq: packing}) hold. Observe that:
\begin{itemize}
\item If $q\in\Dminus$, then since $\chi^2(\mathfrak{q},\mathfrak{u})\le9$ and $\TV(\mathfrak{q},\mathfrak{u})\le\eta$,  
the upper bound in (\ref{eq: packing}) and the decomposition  (\ref{eq: Wp decomp}) imply 
\begin{equation}
\Wp^p(\mu_{F,\mathfrak{q}},\overline{\nu}) \le \underline{\delta}^p + (1-\vartheta)C_0 9^{p/d}\eta^{1-2p/d} M^{-p/d} \le \underline{\delta}^p +  \Delta =: c^p.
\label{eq: lower bound}
\end{equation}
\item On the other hand, if $q\in\Dplus$, then since $\TV(\mathfrak{q},\mathfrak{u})\ge 1/4$, the lower bound in (\ref{eq: packing}) and the decomposition (\ref{eq: Wp decomp}) imply
\begin{equation}
\Wp^p(\mu_{F,\mathfrak{q}},\overline{\nu})\ge\underline{\delta}^p+4\Delta.
\label{eq: upper bound}
\end{equation}
\end{itemize}
Observe that, by Lemma \ref{lem: rootgap},
\[
\begin{split}
s&:=\frac{(\underline{\delta}^p+4\Delta)^{1/p} - (\underline{\delta}^p+\Delta)^{1/p}}{2} \\
&\asymp_{p} \Delta^{1/p} \wedge (\underline{\delta}^{1-p}\Delta) \\
&\asymp_{p,d} (n\log n)^{-1/d} \wedge \big ( \underline{\delta}^{1-p} (n\log n)^{-p/d} \big).
\end{split}
\]

Now, pick any estimator $\hat{W}$ for $\Wp (\cdot,\cdot)$. 
Set $\psi = \mathbbm{1}_{\{\hat{W} \ge c+s\}}$ and view $\psi$ as a randomized test based on $I_1,\dots,I_n \stackrel{\text{i.i.d.}}{\sim} \mathfrak{q}$ (cf. the discussion before Lemma \ref{lem: Wp exact}). Lemma \ref{lem: Wp testing} implies that there exist $\mathfrak{q}^{-} \in \Dminus$ and $\mathfrak{q}^{+} \in \Dplus$ such that
\[
\Prob_{\mathfrak{q}^{-}}(\psi=1) + \Prob_{\mathfrak{q}^{+}}(\psi=0) \ge 0.8.
\]
Recall the notation $\Lambda_{\mathfrak{q}}$ in (\ref{eq: prior}) and use $(\Lambda_0,\Lambda_1) = (\Lambda_{\frakq^-},\Lambda_{\frakq^+})$ as priors, which satisfy
\[
\Lambda_0(\Wp \le c) \ge 0.9 \quad \text{and} \quad \Lambda_1(\Wp \ge c+2s) \ge 0.9,
\]
where we used the fact that each of the events (\ref{eq: lower bound}) and (\ref{eq: upper bound}) holds with probability at least $0.9$ (cf. Lemma \ref{lem: Wp packing}).
The marginal distributions $\overline{P}_0,\overline{P}_1$ satisfy
\[
\begin{split}
&\overline{P}_0 (\hat{W} \ge c+s) + \overline{P}_1 (\hat{W} < c+s)  =\Prob_{\mathfrak{q}^{-}}(\psi=1) + \Prob_{\mathfrak{q}^{+}}(\psi=0) \ge 0.8.
\end{split}
\]
As such, we may apply Proposition \ref{prop: fuzzy} with $\beta_0=\beta_1=0.1$ and $\alpha=0.2$ to conclude that 
\[
\sup_{\Wp(\mu,\nu) \ge \underline{\delta}} \Prob_{\mu,\nu}(|\hat{W}-\Wp(\mu,\nu)| \ge s) \ge 0.3.
\]
This implies the desired result.
\end{proof}

\subsection{Lower bound for $p>2$}
\label{sec: lower bound large p}

The previous construction of random discrete measures yields a lower bound of order $\ell_N \wedge (\underline{\delta}^{1-p} \ell_N^p)$ for $p > 2$, which is loose compared with $\ell_N \wedge (\underline{\delta}^{-1} \ell_N^2)$. Indeed, whenever $\underline{\delta} \ge \ell_N$,
\[
\ell_N \wedge (\underline{\delta}^{1-p} \ell_N^p) = \underline{\delta}^{1-p}\ell_N^p = \underline{\delta} (\ell_N/\underline{\delta})^p \le \underline{\delta} (\ell_N/\underline{\delta})^2 = \ell_N \wedge (\underline{\delta}^{-1} \ell_N^2). 
\]
On the other hand, when $\underline{\delta} < \ell_N$, $\ell_N \wedge (\underline{\delta}^{1-p} \ell_N^p) = \ell_N \wedge (\underline{\delta}^{-1} \ell_N^2) = \ell_N$. 

We consider different random discrete measures, for which the decomposition (\ref{eq: Wp translate-rem}), akin to the one in (\ref{eq: Wp decomp}) in the previous case, holds. The decomposition involves the \textit{Bregman divergence} (cf. Definition 9.2 in \cite{beck2017first}) for $f(z) = \| z \|^p$, 
which  is defined by 
\[
B_f(x,y):=f(x)-f(y)-\langle\nabla f(y),x-y\rangle, \ x,y \in \R^d,
\]
where $\langle \cdot,\cdot \rangle$ denotes the inner product in $\R^d$. A key observation is to construct such random discrete measures to have mean zero, which ensures that the linear term vanishes. Before proceeding, we shall verify the following technical estimate concerning the Bregman divergence. 

\begin{lemma}
\label{lem: bregman}
For all $z$ and $v$ in $\mathbb R^d$,
\begin{equation}
 B_f(z+v,z)\asymp_p \|z\|^{p-2}\|v\|^2+\|v\|^p.
 \label{eq: bregman}
\end{equation}
\end{lemma}

\begin{proof}
This follows from equation (2.19) in \cite{pinchover2008a}. For completeness, we provide a proof.
Since $p>2$, the function $f(x)=\|x\|^p$ is twice continuously differentiable with $\nabla^2 f(0)=0$. For $x\ne0$,
\[
 \langle v,\nabla^2 f(x)v\rangle
 =p\|x\|^{p-2}\|v\|^2+p(p-2)\|x\|^{p-4}\langle x,v\rangle^2.
\]
By the Cauchy--Schwarz inequality, $\langle v,\nabla^2 f(x)v\rangle$ lies between $p\|x\|^{p-2}\|v\|^2$ and $p(p-1)\|x\|^{p-2}\|v\|^2$. Both bounds also hold at $x=0$. Applying Taylor's formula with integral remainder to $t\mapsto f(z+tv)$ between $0$ and $1$, we obtain
\[
 \begin{split}
 B_f(z+v,z)
 &=f(z+v)-f(z)-\langle\nabla f(z),v\rangle\\
 &=\int_0^1(1-t)\langle v,\nabla^2 f(z+tv)v\rangle\,dt\\
 &\asymp_p \|v\|^2\int_0^1(1-t)\|z+tv\|^{p-2}\,dt.
 \end{split}
\]
Since $\|z+tv\|\le\|z\|+\|v\|$ for $0\le t\le1$, the last integral is at most $(\|z\|+\|v\|)^{p-2}/2$.

For the lower bound, suppose first that $\|v\|\le2\|z\|$. For $0\le t\le1/4$, the reverse triangle inequality gives
\[
 \|z+tv\| \ge \|z\|-t\|v\| \ge \frac{\|z\|}{2} \ge \frac{\|z\|+\|v\|}{6}.
\]
If instead $\|v\| > 2\|z\|$, then for $3/4 \le t \le1$,
\[
 \|z+tv\| \ge t\|v\|-\|z\| \ge \frac{\|v\|}{4} \ge \frac{\|z\|+\|v\|}{6}.
\]
In either case, we have
\[
 \int_0^1(1-t)\|z+tv\|^{p-2}\,dt
 \gtrsim \left(\frac{\|z\|+\|v\|}{6}\right)^{p-2}.
\]
Combining the two bounds and using $(\|z\|+\|v\|)^{p-2}\asymp_p\|z\|^{p-2}+\|v\|^{p-2}$, we obtain (\ref{eq: bregman}).
\end{proof}

We are ready to prove the lower bound for $p > 2$.

\begin{proof}[Proof of Theorem \ref{thm: Wp main}: Lower bound for $p > 2$]
As before, we assume $N = n \le m$. 
Let $M \ge 2$ be a positive integer to be determined later. 
 Recall that $\mathfrak{u}$ is the uniform distribution on $[M]$. Let $F:[M] \to \calS$ be a uniform random bijection from Lemma \ref{lem: Wp packing} with $\calB=B_d(e_1/6,1/12)$.

\medskip

\underline{Step 1} (Construction of random discrete measures). 
For any $\mathfrak{q} \in \calP([M])$, consider
\[
 \mu_{F,\mathfrak{q}}=\sum_{i=1}^M \frakq_i \frac{\delta_{F(i)}+\delta_{-F(i)}}2,
\]
where $\frakq_i = \frakq (\{i\})$.
Define $\overline{\nu}$ by taking $\mathfrak{q}=\mathfrak{u}$,
\[
\overline{\nu} = \frac{1}{M} \sum_{i=1}^M \frac{\delta_{F(i)}+\delta_{-F(i)}}2 = \frac{1}{2M} \sum_{x \in \calS \cup (-\calS)} \delta_{x},
\]
which does not depend on $F$. 
Let $z=\underline{\delta} e_1$ and $\mathfrak{t}_z(y)=y+z$. Since $\|F(i)\|\le 1/4$, both $\mu_{F,\mathfrak{q}}$ and $(\mathfrak{t}_{z})_{\#}\overline{\nu}$ are supported in $B_d(0,1)$ for $\underline{\delta} \le 3/4$.

For any $\pi \in \Pi(\mu_{F,\mathfrak{q}},\overline{\nu})$, symmetry of the two marginals gives $\int(y-x)\,d\pi(x,y)=0$. Integrating the identity 
\[
f(z+y-x)=f(z)+\langle\nabla f(z),y-x\rangle+B_f(z+y-x,z)
\]
with respect to $\pi$ and taking the infimum over $\pi$, we obtain
\begin{equation}
 T_{\mathfrak{q}}(F):=\Wp^p(\mu_{F,\mathfrak{q}},(\mathfrak{t}_{z})_{\#}\overline{\nu})
 =\underline{\delta}^p+
    \inf_{\pi\in\Pi(\mu_{F,\mathfrak{q}},\overline{\nu})}\int B_f(z+y-x,z)\,d\pi(x,y),
    \label{eq: Wp translate-rem}
\end{equation}
where we used the fact that the pushforward under $(x,y)\mapsto(x,y+z)$ is a bijection from $\Pi(\mu_{F,\mathfrak{q}},\overline{\nu})$ onto $\Pi(\mu_{F,\mathfrak{q}},(\mathfrak{t}_{z})_{\#}\overline{\nu})$.
By convexity of $f$, $B_f(z+y-x,z)\ge0$, which guarantees that $T_{\mathfrak{q}}(F)\ge\underline{\delta}^p$. Set
\[
\Delta=\underline{\delta}^{p-2}M^{-2/d}+M^{-p/d}.
\]

\medskip

\underline{Step 2} (Lower bound for $ T_{\mathfrak{q}}(F)$ under $\Dplus$). 
Let $\mathcal C_i=\{-F(i),F(i)\}$ for $i\in[M]$. We first bound the cost of moving mass between distinct sets $\mathcal C_i$ and $\mathcal C_j$. For $i\ne j$, Lemma \ref{lem: Wp packing} gives
$\|F(i)-F(j)\|\ge c_0 M^{-1/d}$. Since $F(i)$ and $F(j)$ lie in
$B_d(e_1/6,1/12)$, we also have $\|F(i)+F(j)\|\ge 2(1/6-1/12)=1/6$. Since $M^{-1/d} \le1$,
\[
\dist(\mathcal{C}_i,\mathcal{C}_j) \gtrsim_d M^{-1/d}.
\]
Any $\pi\in\Pi(\mu_{F,\mathfrak{q}},\overline{\nu})$ satisfies
$\pi(\mathcal{C}_i\times \mathcal{C}_i) \le \min \{ \mathfrak{q}_i,1/M \}$. For $q\in\Dplus$, the union of
$\mathcal{C}_i\times \mathcal{C}_j$ with $i\ne j$ therefore has $\pi$-mass at least
\[
1-\sum_{i \in [M]} \pi(\calC_i \times \calC_i) \ge 1-\sum_{i \in [M]} \min \{ \mathfrak{q}_i,1/M \} = \TV(\mathfrak{q},\mathfrak{u}) \ge \frac{1}{4}.
\]
The penultimate equality follows by Scheff\'{e}'s lemma (cf. (\ref{eq: scheffe})).
For $(x,y) \in \bigcup_{i \ne j} \mathcal{C}_i\times \mathcal{C}_j$, the estimate in (\ref{eq: bregman}) yields 
\[
\begin{split}
B_f(z+y-x,z) &\gtrsim_{p} \|z\|^{p-2}\|y-x\|^2+\|y-x\|^p \\
&\gtrsim_d \underline{\delta}^{p-2}M^{-2/d} + M^{-p/d} = \Delta.
\end{split}
\]
As such, restricting the integral in 
(\ref{eq: Wp translate-rem}) to $\bigcup_{i \ne j} \mathcal{C}_i\times \mathcal{C}_j$, we obtain
\[
T_{\mathfrak{q}}(F)\ge\underline{\delta}^p+c_2\Delta
\]
with $c_2 \gtrsim_{p,d} 1$.

\medskip

\underline{Step 3} (Upper bound for $T_{\mathfrak{q}}(F)$ under $\Dminus$).
For any fixed $q\in\Dminus$, by Lemma \ref{lem: Wp packing}, it holds that
\begin{equation}
\Wp^p(F_{\#}\mathfrak{q},F_{\#}\mathfrak{u})\le L_\eta M^{-p/d}
\label{eq: Wass event}
\end{equation}
with probability at least $0.9$, where $L_\eta=10C_09^{p/d}\eta^{1-2p/d}$. Given $F$, generate $(U,V)$ from an optimal coupling for $\Wp(F_{\#}\mathfrak{q},F_{\#}\mathfrak{u})$. Let $\varepsilon$ be an independent Rademacher random variable, i.e., $\Prob(\varepsilon = \pm 1)=1/2$. Then $(\varepsilon U,\varepsilon V)$ has marginals $\mu_{F,\mathfrak{q}}$ and $\overline{\nu}$, and $\|\varepsilon U-\varepsilon V\|=\|U-V\|$. Condition on the event (\ref{eq: Wass event}). We have
\[
\begin{split}
&\E[\|U-V\|^p\mid F]=\Wp^p(F_{\#}\mathfrak{q},F_{\#}\mathfrak{u})\le L_\eta M^{-p/d}, \\
&\E[\|U-V\|^2\mid F]\le L_\eta^{2/p}M^{-2/d},
\end{split}
\]
where the second inequality follows by H\"older's inequality. Using the conditional distribution of $(\varepsilon U,\varepsilon V)$ in
(\ref{eq: Wp translate-rem}), combined with (\ref{eq: bregman}), we obtain
\[
 T_{\mathfrak{q}}(F)-\underline{\delta}^p \le \E\big [  B_f(z+\varepsilon (V-U),z) \mid F\big ] 
 \le C_2 \big ( \underline{\delta}^{p-2}L_\eta^{2/p}M^{-2/d}+L_\eta M^{-p/d} \big),
\]
where $C_2 \lesssim_p 1$. 
Since $1-2p/d>0$,  one may choose 
sufficiently small $0<\eta\le1/10$ (that may depend on $(p,d)$) such that 
\[
C_2\max\{ L_\eta^{2/p},L_\eta \}\le \frac{c_2}{2},
\]
where $c_2$ is given in Step 2. Then,
\[
\Prob_F\big(T_{\mathfrak{q}}(F)\le\underline{\delta}^p+c_2\Delta/2\big) \ge 0.9
\]
for any
$q\in\Dminus$. 

Choose $M$ as in (\ref{eq: M-choice}) with $k=n$.
Given that 
\[
\begin{split}
(\underline{\delta}^p + c_2\Delta)^{1/p} - (\underline{\delta}^p + c_2\Delta/2)^{1/p} 
&\asymp_{p,d} \Delta^{1/p} \wedge (\underline{\delta}^{1-p}\Delta) \\
&\gtrsim M^{-1/d} \wedge (\underline{\delta}^{-1}M^{-2/d}) \\
&\asymp_{p,d} (n\log n)^{-1/d} \wedge \big ( \underline{\delta}^{-1} (n\log n)^{-2/d} \big)
\end{split}
\]
the rest of the proof is similar to the $1 \le p \le 2$ case. 
\end{proof}

\subsection{Upper bound}
\label{sec: Wp upper}

The upper bound for $p > 2$ rests on the following proposition, whose proof is deferred to the next section.

\begin{proposition}
\label{prop: power bound}
Let $d > 2p > 4$. Pick any $(\mu,\nu) \in \calP(B_d(0,1)) \times \calP(B_d(0,1))$. For $t=\Wp(\mu,\nu)$, we have 
\begin{equation}
\E[|\Wp^p(\hat{\mu}_n,\hat{\nu}_m)-t^p|] \lesssim_{p,d} t^{p-2} N^{-2/d} + N^{-p/d}.
\label{eq: power bound}
\end{equation}
\end{proposition}

Given this proposition, the upper-bound proof is straightforward. 

\begin{proof}[Proof of Theorem \ref{thm: Wp main}: Upper bound]
Let $t=\Wp(\mu,\nu) \ge \underline{\delta}$ and $\hat{W}=\Wp(\hat\mu_n,\hat\nu_m)$. By symmetry, we may assume $N=n \le m$. 
We shall establish $\E[|\hat{W}-t|] \lesssim_{p,d} n^{-1/d}\wedge(\underline{\delta}^{1-\alpha}n^{-\alpha/d})$. By the triangle inequality, Jensen's inequality, and Theorem 1 in \cite{fournier2015on}, we have
\[
 \E[|\hat{W}-t|]
 \le \E[\Wp(\hat\mu_n,\mu)]+\E[\Wp(\hat\nu_m,\nu)]
 \lesssim_{p,d} n^{-1/d}+m^{-1/d} \lesssim n^{-1/d}.
\] 
As such, it remains to establish $\E[|\hat{W}-t|] \lesssim_{p,d} \underline{\delta}^{1-\alpha}n^{-\alpha/d}$ when $t \ge \underline{\delta} \ge n^{-1/d}$. 
For $1 \le p \le 2$, the result is known; see Corollary 4 and Remark 1 in \cite{manole2024sharp}. See also Remark \ref{rem: upper bound for small p} below.  We shall focus on the $p > 2$ case. 
We shall use Proposition \ref{prop: power bound}. 
When  $t\ge\underline{\delta}\ge n^{-1/d}$, $t^{p-2}n^{-2/d}+n^{-p/d} \lesssim t^{p-2}n^{-2/d}$. We conclude that
\[
\E[|\hat{W}-t|]\le t^{1-p}\E[|\hat{W}^p-t^p|] \lesssim_{p,d} t^{-1} n^{-2/d} \le \underline{\delta}^{-1}n^{-2/d}.
\]
as desired. 
\end{proof}

\subsection{Proof of Proposition \ref{prop: power bound}}
\label{sec: power bound}
As before, we assume $N = n \le m$. The proof is divided into three steps. 

\medskip
\underline{Step 1} (Reduction to equal sample size).
Let $\hat{T}_{k,\ell}=\Wp^p(\hat\mu_k,\hat\nu_\ell)$, using the first $k$ observations from $\mu$ and the first $\ell$ from $\nu$. We shall verify that
\[
t^p\le\E[\hat{T}_{n,m}] \le \E[\hat{T}_{n,n}].
\]
Pick any continuous functions $f$ and $g$ on $\calX$ satisfying $f(x)+g(y)\le\|x-y\|^p$ for all $x$ and $y$ in $\calX$. By duality (cf. (\ref{eq: OT-dual})), we have
\[
 \int f\,d\mu+\int g\,d\nu
 =\E\left[ \int f \, d\hat{\mu}_n + \int g \, d\hat{\nu}_m \right] 
 \le \E[\hat{T}_{n,m}].
\]
Taking the supremum over all such $f$ and $g$ on the left-hand side yields $t^p\le\E[\hat{T}_{n,m}]$.

Let $\mathcal J$ be the set of $n$-element subsets of $[m]$ and put $\hat\nu_J=n^{-1}\sum_{j\in J}\delta_{Y_j}$ for $J\in\mathcal J$. Each index $j$ belongs to $\binom{m-1}{n-1}$ elements of $\mathcal J$. Hence,  
\[
 \hat\nu_m=\frac{1}{|\mathcal J|}\sum_{J\in\mathcal J}\hat\nu_J \quad \text{and} \quad
 \hat{T}_{n,m}\le\frac{1}{|\mathcal J|}\sum_{J\in\mathcal J}
      \Wp^p(\hat\mu_n,\hat\nu_J),
\]
where the inequality follows from convexity of $\Wp^p$ in the second marginal. For each fixed $J$, $\Wp^p(\hat\mu_n,\hat\nu_J) \stackrel{d}{=} \hat{T}_{n,n}$, where $\stackrel{d}{=}$ denotes equality in distribution. Taking expectation, we obtain  $\E[\hat{T}_{n,m}]\le\E[\hat{T}_{n,n}]$, as desired.

\medskip
\underline{Step 2} (Bounding variance). 
Let $(X'_1,\ldots,X'_n,Y'_1,\ldots,Y'_m)$ be an independent copy of the data. Pick any $i \in [n]$. Let $\hat{\mu}_n'$ be the empirical distribution of $\{ X_j \}_{j \in [n] \setminus \{ i \}} \cup \{ X_i' \}$. Using the expression (\ref{eq: discrete Wp}), we observe that
\[
\big|\Wp^p(\hat{\mu}_n,\hat{\nu}_m) - \Wp^p(\hat{\mu}_n',\hat{\nu}_m)\big| \le \sup_{\gamma} 
 \left|\sum_{j=1}^m\gamma_{ij}
    \bigl(\|X_i'-Y_j\|^p-\|X_i-Y_j\|^p\bigr)\right|
 \le \frac{2^p}{n}.
\]
Likewise, for any $j \in [m]$, let $\hat{\nu}_m'$ be the empirical distribution of $\{ Y_i \}_{i \in [m] \setminus \{ j \}} \cup \{ Y_j' \}$; then 
\[
\big|\Wp^p(\hat{\mu}_n,\hat{\nu}_m) - \Wp^p(\hat{\mu}_n,\hat{\nu}_m')\big| 
 \le \frac{2^p}{m}.
\]
Now, the Efron-Stein inequality (cf. Theorem 3.1 in \cite{boucheron2013concentration}) yields
\[
\Var(\hat{T}_{n,m})
 \le \frac{2^{2p}}{2}\left(\frac{1}{n}+\frac{1}{m}\right)
 \le \frac{2^{2p}}{n}.
\]
Since $\E[|\hat{T}_{n,m}-\E[\hat{T}_{n,m}]|]\le\sqrt{\Var(\hat{T}_{n,m})}$, it follows that
\begin{equation}
 \E[|\hat{T}_{n,m}-t^p|]\le \E[|\hat{T}_{n,m} - \E[\hat{T}_{n,m}]|] + \E[\hat{T}_{n,m}] - t^p \le 2^p n^{-1/2}+\E[\hat{T}_{n,n}]-t^p.
 \label{eq: Wp cost-bias}
\end{equation}

\medskip

\underline{Step 3} (Bounding bias). Recall that we have assumed $p > 2$. Let $f(z) = \|z\|^p$. 
We shall establish
\begin{equation}
\E[\hat{T}_{n,n}]-t^p\lesssim_{p,d}t^{p-2}n^{-2/d}+n^{-p/d}.
 \label{eq: Wp local-cost}
\end{equation}
Let $\pi$ be an optimal coupling for $\Wp(\mu,\nu)$. Generate two independent samples $\{ (U_i,V_i) \}_{i=1}^n$ and $\{ (U'_j,V'_j) \}_{j=1}^n$ from $\pi$. We shall bound $\E[\hat{T}_{n,n}]$ by matching the first coordinates of each sample to the second coordinates of the other.

Let $Z_i=(U_i+V_i)/2$ and $Z_j'=(U_j'+V_j')/2$. By (\ref{eq: matching}), we have
\begin{equation}
 S:=\Wp^p \left(\frac{1}{n}\sum_{i=1}^n\delta_{Z_i},
                  \frac{1}{n}\sum_{j=1}^n\delta_{Z_j'}\right) = \min_{\sigma}\frac{1}{n}\sum_{i=1}^n\|Z_i-Z'_{\sigma(i)}\|^p, 
\label{eq: permutation}
\end{equation}
where the minimum is taken over all permutations $\sigma$ of $[n]$. By the triangle inequality and Theorem 1 in \cite{fournier2015on}, we have $\E[S] \lesssim_{p,d} n^{-p/d}$. 

Let $\sigma$ be any minimizer in (\ref{eq: permutation}), which depends on $(Z_1,\dots,Z_n)$ and $(Z_1',\dots,Z_n')$. 
Let 
\[
\Delta_i=U_i-V_i \quad \text{and} \quad \Delta_j'=U_j'-V_j'.
\]
For $i\in[n]$, put 
\[
a_i=(\Delta_i+\Delta'_{\sigma(i)})/2 \quad \text{and} \quad b_i=Z_i-Z'_{\sigma(i)}.\]
Then $U_i-V'_{\sigma(i)}=a_i+b_i$ and $U'_{\sigma(i)}-V_i=a_i-b_i$.
Assigning mass $1/n$ to each $(U_i,V'_{\sigma(i)})$ couples the empirical distributions of $\{ U_i \}_{i=1}^n$ and $\{ V_j' \}_{j=1}^n$, so that 
\[
\hat{T}_{n,n} \stackrel{d}{=} \Wp^p \left(\frac{1}{n}\sum_{i=1}^n\delta_{U_i},
                  \frac{1}{n}\sum_{j=1}^n\delta_{V_j'}\right)  \le \frac{1}{n}\sum_{i=1}^n f(a_i+b_i).
\]
Likewise, assigning mass $1/n$ to each $(U_{\sigma(i)}',V_i)$ couples the empirical distributions of  $\{ U_j' \}_{j=1}^n$ and $\{ V_i \}_{i=1}^n$, so that 
\[
\hat{T}_{n,n} \stackrel{d}{=} \Wp^p \left(\frac{1}{n}\sum_{i=1}^n\delta_{U_i'},
                  \frac{1}{n}\sum_{j=1}^n\delta_{V_j}\right)  \le \frac{1}{n}\sum_{i=1}^n f(a_i-b_i).
\]
As such, we obtain
\begin{equation}
\E[\hat{T}_{n,n}]\le\E\left[\frac{1}{2n}\sum_{i=1}^n \big ( f(a_i+b_i) + f(a_i-b_i) \big ) \right].
\label{eq: bias}
\end{equation}

Recall the Bregman divergence: $B_f(x,y) = f(x)-f(y) - \langle \nabla f(y),x-y \rangle$. 
Observe that
\[
 f(a_i+b_i) + f(a_i-b_i)=2f(a_i)+B_f(a_i+b_i,a_i)+B_f(a_i-b_i,a_i).
\]
Convexity of $f$ yields $2f(a_i)\le f(\Delta_i)+f(\Delta'_{\sigma(i)})$. Combining (\ref{eq: bias}), we have 
\[
\begin{split}
\E[\hat{T}_{n,n}] &\le \E\left [ \frac{1}{2n} \sum_{i=1}^n \big (f(\Delta_i)+f(\Delta_{\sigma(i)}') \big )  \right ] + \E\left [ \frac{1}{2n} \sum_{i=1}^n \big (B_f(a_i+b_i,a_i)+B_f(a_i-b_i,a_i) \big )  \right ] \\
&=\underbrace{\frac{1}{n}\sum_{i=1}^n \E[f(\Delta_i)]}_{=t^p} + \E\left [ \frac{1}{2n} \sum_{i=1}^n \big (B_f(a_i+b_i,a_i)+B_f(a_i-b_i,a_i) \big )  \right ],
\end{split}
\]
where we used the fact that $\E[f(\Delta_i)] = t^p$, which follows by construction. 

By Lemma \ref{lem: bregman}, 
\[
B_f(a_i+b_i,a_i)+B_f(a_i-b_i,a_i) \lesssim_p \|a_i\|^{p-2}\|b_i\|^2+\|b_i\|^p.
\]
Since $S=n^{-1}\sum_{i=1}^n\|b_i\|^p$, H\"older's inequality  yields
\[
 \E\left[\frac{1}{n}\sum_{i=1}^n\|a_i\|^{p-2}\|b_i\|^2\right]
 \le \left(\E\left [\frac{1}{n}\sum_{i=1}^n f(a_i) \right]\right)^{(p-2)/p}
       \left(\E[S]\right)^{2/p}
 \le t^{p-2}\left(\E[S]\right)^{2/p}.
\]
Using $\E[S]\lesssim_{p,d} n^{-p/d}$, we conclude
\[
 \E[\hat{T}_{n,n}]-t^p
 \lesssim_p t^{p-2}(\E[S])^{2/p}+\E[S]
 \lesssim_{p,d}t^{p-2}n^{-2/d}+n^{-p/d}.
\]
Finally, combining (\ref{eq: Wp local-cost}) and (\ref{eq: Wp cost-bias}), we obtain (\ref{eq: power bound}). This completes the proof. 
\qed

\begin{remark}
\label{rem: upper bound for small p}
The same argument gives the bound for $1\le p\le2$ for $d>2p$. Indeed, the construction leading to (\ref{eq: bias}) applies for all $p\ge1$. By concavity of $u\mapsto u^{p/2}$ and the parallelogram identity,
\[
 \frac{\|a+b\|^p+\|a-b\|^p}{2}
 \le (\|a\|^2+\|b\|^2)^{p/2}
 \le \|a\|^p+\|b\|^p.
\]
Taking $a=a_i$ and $b=b_i$ in (\ref{eq: bias}) and using convexity of $\|\cdot\|^p$, we have $\E[\hat{T}_{n,n}]-t^p\le\E[S]\lesssim_{p,d} n^{-p/d}$. Since $n^{-1/2}\le n^{-p/d}$ for $d>2p$, (\ref{eq: Wp cost-bias}) yields $\E[|\hat W^p-t^p|]\lesssim_{p,d} n^{-p/d}$. Since $|\hat W-t|\le t^{1-p}|\hat W^p-t^p|$, the desired bound follows.
\end{remark}

\section{Proof of Theorem \ref{thm: Wp main small p}}
\label{sec: proof Wp small p}

\subsection{Preliminaries}
\label{sec: Wp small p prelim}

The lower-bound proof for Theorem \ref{thm: Wp main} used the $d > 2p$ condition when it invoked  Lemma \ref{lem: Wp packing} (Proposition 9 in \cite{nilesweed2022estimation}), more specifically, the upper bound in (\ref{eq: packing}). The proof of Proposition 9 in \cite{nilesweed2022estimation} rests on bounding $\Wp$ using dyadic partitions; cf. \cite{fournier2015on, weed2019sharp}. Its proof can be extended to $d \le 2p$, but the lower-bound proof for $4 < d \le 2p$ requires a sharper treatment. Recalling the remainder estimate (\ref{eq: bregman}) for the Bregman divergence, one will need to suitably upper bound the transport cost for a \textit{mixed-power} ground cost of the form $\mathsf{c} (x,y) = b\|x-y\|^2 + \|x-y\|^p$ for some $b>0$ and keep track of the dependence on $b$. The following estimate, which builds on dyadic partitions and is adapted from Lemma 3.3 in \cite{fournier2023convergence}, will suffice for our purpose. Lemma \ref{lem: dyadic} will also play a key role in the upper-bound proof.

\begin{lemma}[Dyadic bound for general cost]
\label{lem: dyadic}
Let $\calQ_\ell$ be the dyadic partition of $\mathsf{D}=[-1,1]^d$ into cubes of
side length $2^{1-\ell}$. Let $\sfc: \mathsf{D} \times \mathsf{D} \to \R_+$ be a continuous cost with $\sfc(x,x) = 0$ for $x \in \mathsf{D}$. 
Pick any $\mu,\nu \in\calP(\mathsf{D})$. For $Q \in \calQ_\ell$, put $\overline{\sfc}_Q=\sup_{x,y\in Q}\sfc(x,y), a_Q=|\mu(Q)-\nu(Q)|,  \overline{\sfC}_\ell=\max_{Q\in\calQ_\ell}\overline{\sfc}_Q$, and $A_\ell=\sum_{Q \in \calQ_\ell}a_Q$. 
Denote $Q'\prec Q$ if $Q' \in \calQ_{\ell+1}$ and $Q'\subset Q$.
Then there exists $\pi\in\Pi(\mu,\nu)$, independent of $\sfc$, such that
\begin{equation}
 \int\sfc\,d\pi
 \le \sum_{\ell=0}^{\infty}\sum_{Q\in\calQ_\ell}
      \overline{\sfc}_Q\sum_{Q'\prec Q}a_{Q'}
 \le \sum_{\ell=0}^{\infty}\overline{\sfC}_\ell A_{\ell+1}.
 \label{eq: dyadic moment bound}
\end{equation}
\end{lemma}

\begin{remark}
For $\sfc(x,y)=\|x-y\|^r$ for some $r > 0$, $\overline{\sfC}_\ell \lesssim_{r,d} 2^{-r\ell}$. 
In our applications, we combine Lemma \ref{lem: dyadic} with the following estimate on a dyadic series: for $a,b,r>0$ with $a/b \ge 1$, 
\begin{equation}
 \sum_{\ell\ge1}2^{-r\ell}\big(a \wedge b\,2^{d\ell/2}\big)
 \lesssim_{r,d}
 \begin{cases}
 b, & d<2r,\\
 b\bigl(1+\log(a/b)\bigr), & d=2r,\\
 a^{1-2r/d}b^{2r/d}, & d>2r.
 \end{cases}
 \label{eq: dyadic series}
\end{equation}
This follows by summing separately over $2^{d\ell/2}\le a/b$ and $2^{d\ell/2}>a/b$.
\end{remark}

\begin{proof}[Proof of Lemma \ref{lem: dyadic}]
The proof is a mild modification of that of Lemma 3.3 in \cite{fournier2023convergence}. We provide an outline of the proof.

\medskip

\underline{Step 1} (Coupling construction). 
For dyadic partitions $\calQ_0,\dots,\calQ_k$, let $\pi_k$ be a coupling for $(\mu,\nu)$ constructed in Steps 1-2 in the proof of Lemma 3.3 in \cite{fournier2023convergence}. We review its construction below. For $Q\in\calQ_k$ with $\mu(Q)\nu(Q)>0$, put
\[
\pi_Q
=\frac{\mu|_Q \otimes \nu|_Q}{\mu(Q)\nu(Q)}.
\]
 Proceeding by reverse induction, let $Q \in \calQ_\ell$,
$\ell<k$, with $\mu(Q)\nu(Q)>0$.  For each $Q' \prec Q$, define
\[
\rho_{Q'}=\frac{\mu(Q')}{\mu(Q)}
                  \wedge\frac{\nu(Q')}{\nu(Q)},
\quad
q_Q=1-\sum_{Q'\prec Q}\rho_{Q'}
=\frac{1}{2}\sum_{Q'\prec Q}
\left|\frac{\mu(Q')}{\mu(Q)}-
            \frac{\nu(Q')}{\nu(Q)}\right|.
\]
If $q_Q>0$, set
\begin{align*}
\alpha_Q
&=\frac1{q_Q}\sum_{Q'\prec Q}
 \left(\frac{\mu(Q')}{\mu(Q)}-
             \frac{\nu(Q')}{\nu(Q)}\right)_+
 \frac{\mu|_{Q'}}{\mu(Q')},\\
\beta_Q
&=\frac1{q_Q}\sum_{Q'\prec Q}
 \left(\frac{\nu(Q')}{\nu(Q)}-
             \frac{\mu(Q')}{\mu(Q)}\right)_+
 \frac{\nu|_{Q'}}{\nu(Q')}.
\end{align*}
These measures are supported in $Q$. Define
\begin{equation}\label{eq:recursion}
\pi_Q
=\sum_{Q'\prec Q}\rho_{Q'}\pi_{Q'}
  +q_Q\alpha_Q \otimes \beta_Q.
\end{equation}
Discard the last term when $q_Q=0$. 

By induction, the first marginal of \eqref{eq:recursion} is
\[
\sum_{Q'\prec Q}
\left[\rho_{Q'}+
 \left(\frac{\mu(Q')}{\mu(Q)}-
             \frac{\nu(Q')}{\nu(Q)}\right)_+\right]
\frac{\mu|_{Q'}}{\mu(Q')}
=\frac{\mu|_Q}{\mu(Q)}.
\]
Likewise, the second marginal is $\nu|_Q/\nu(Q)$. Let $\pi_k = \pi_{\mathsf{D}} \in \Pi(\mu,\nu)$. The construction is independent of $\sfc$.

\medskip

\underline{Step 2} (Moment bound).
Set $u_{\mathsf{D}}=1$ and $u_{Q'}=u_Q\rho_{Q'}$ for $Q'\prec Q$, with
$u_{Q'}=0$ if $u_Q=0$. By the definition of $\rho_{Q'}$ and induction, we have
$u_Q\le\mu(Q)\wedge\nu(Q)$. Hence, for $u_Q>0$,
\[
\begin{split}
 2u_Qq_Q
 &\le \sum_{Q'\prec Q}\left|\mu(Q')-\frac{\mu(Q)}{\nu(Q)}\nu(Q')\right|\\
 &\le \sum_{Q'\prec Q}a_{Q'}+a_Q
 \le 2\sum_{Q'\prec Q}a_{Q'},
\end{split}
\]
where we used $\sum_{Q'\prec Q}\nu(Q')=\nu(Q)$.

Iterating \eqref{eq:recursion}, we obtain
\[
\begin{split}
 \int\sfc\,d\pi_k
 &\le \sum_{Q\in\calQ_k}u_Q\overline{\sfc}_Q
      +\sum_{\ell=0}^{k-1}\sum_{Q\in\calQ_\ell}u_Qq_Q\overline{\sfc}_Q\\
 &\le \overline{\sfC}_k
      +\sum_{\ell=0}^{k-1}\sum_{Q\in\calQ_\ell}
       \overline{\sfc}_Q\sum_{Q'\prec Q}a_{Q'},
\end{split}
\]
where we used $\sum_{Q\in\calQ_k}u_Q \le 1$. 

Finally, by weak compactness of $\Pi(\mu,\nu)$, choose
a subsequence of $\pi_k$ with weak limit $\pi$.
Since $\sfc$ is bounded and continuous and $\overline{\sfC}_k\to 0$, passing
to the limit gives the first inequality in \eqref{eq: dyadic moment bound}.
The second follows from $\overline{\sfc}_Q\le\overline{\sfC}_\ell$
and $\sum_{Q\in\calQ_\ell}\sum_{Q'\prec Q}a_{Q'}=A_{\ell+1}$.
\end{proof}

Another ingredient of our upper-bound proof for $4 < d \le 2p$ is the following technical lemma, which is a consequence of the \textit{shortening} principle for OT; see Chapter 8 in \cite{villani2009optimal}. 

\begin{lemma}[A consequence of shortening principle]
\label{lem: shortening}
Let $p>1$.
Let $\mu,\nu$ be compactly supported, 
and let $\pi$ be an optimal coupling for $\Wp$ with support
$\Gamma:=\spt(\pi)$. Set $\mathcal{Z}=\left\{(u+v)/2:(u,v)\in\Gamma\right\}$. 
Then there exists a Lipschitz function $D:\mathcal{Z}\to\R^d$
such that
\[
D\left(\frac{u+v}{2}\right)=u-v, \quad (u,v)\in\Gamma.
\]
One can extend $D$ to $\R^d$ while keeping the same Lipschitz constant. For $p \ge 2$, $D$ is $2\sqrt{p-1}$-Lipschitz. 
\end{lemma}
\begin{proof}
By optimality, $\Gamma$ is $\sfc$-cyclically monotone with $\sfc(x,y)=\|x-y\|^p$. 
For $1 < p < 2$, the conclusion follows from Theorem 8.23 in
\cite{villani2009optimal}. Its proof extends to $p > 2$.  We provide a short proof of the lemma for $p \ge 2$ with an explicit Lipschitz constant.   

For $a,b\in\R^d$, Jensen's inequality and Theorem 1 in \cite{ball1994sharp} yield
\[
\|a\|^2+\|b\|^2
\le Q_p(a,b):=
\left(\frac{\|a+b\|^p+\|a-b\|^p}{2}\right)^{2/p}
\le \|a\|^2+(p-1)\|b\|^2.
\]
Take $(u,v),(u',v')\in\Gamma$ and set $z=(u+v)/2, z' = (u'+v')/2, a=((u-v)+(u'-v'))/2$, and $w=((u-v)-(u'-v'))/2$.
Cyclical monotonicity of $\Gamma$
implies $Q_p(a,w)\le Q_p(a,z-z')$. As such, 
\[
\|a\|^2+\|w\|^2
\le Q_p(a,w)
\le Q_p(a,z-z')
\le \|a\|^2+(p-1)\|z-z'\|^2,
\]
which implies
\[
\|(u-v)-(u'-v')\| \le 2\sqrt{p-1} \|z-z'\|.
\]
The final claim follows by the Kirszbraun extension theorem (cf. Theorem 2.10.43 in \cite{federer2014geometric}). 
\end{proof}

\begin{remark}[Connection with the shortening principle]
To see how Lemma \ref{lem: shortening} follows from the shortening principle, pick any $(u_1,v_1), (u_2,v_2) \in \Gamma$, and consider their interpolations 
\[
\gamma_j(t) = (1-t)u_j + tv_j, \ t \in [0,1],  j=1,2.
\]
The shortening principle in Chapter 8 in \cite{villani2009optimal} states that one can control the uniform distance between the \textit{whole paths} $\gamma_1, \gamma_2$ by their distance at \textit{any} $t_0 \in (0,1)$, that is, for some constant $C$ (that may depend on $t_0$),
\[
\sup_{t \in [0,1]} \| \gamma_1(t)-\gamma_2(t)\| \le C \| \gamma_1(t_0) - \gamma_2(t_0)\|. 
\]
Taking $t_0=1/2$ implies 
\[
\begin{split}
\|(u_1,v_1) - (u_2,v_2)\| &\le \| u_1-u_2 \|+\|v_1-v_2\| \\
&\le 2C\| \gamma_1(1/2) - \gamma_2(1/2)\|.
\end{split}
\]
\end{remark}

Our upper-bound proof for $4 < d \le 2p$ will adapt some arguments in the proof of Proposition \ref{prop: power bound}. Recalling the notation used in Step 3 of that proof, Lemma \ref{lem: shortening} allows us to represent $\Delta_i$ as $\Delta_i = D(Z_i)$, and the problem reduces to bounding the empirical transport cost for $\sfc(z,z') = (\| D(z) \|  \vee \|D(z')\|)^{p-2}\|z-z'\|^2 + \|z-z'\|^p$. The final step will use Lemma \ref{lem: dyadic}. 

\subsection{Lower bound}

The lower bound follows by combining Propositions \ref{prop: Wp two point} and  \ref{prop: Wp appendix high lower} below. 
\begin{proposition}[Preliminary lower bound]
\label{prop: Wp two point}
Under $n \le m$, for any $0 < \underline{\delta} \le 1$, we have
\begin{equation}
 \calR^{\Wp}_{n,m,d}(\underline{\delta})
 \gtrsim_p n^{-1/(2p)}\wedge(\underline{\delta}^{1-p}n^{-1/2}).
 \label{eq: Wp two-point}
\end{equation}
\end{proposition}
\begin{proof}
By Lemma \ref{lem: minimax risk}~(iii), 
it suffices to consider $d=1$. Let $\nu=(1/4)\delta_1+(3/4)\delta_{-1}$, and put
\[
\mu_t = \left(\frac14+t\right)\delta_1+\left(\frac34-t\right)\delta_{-1}, \quad 0\le t\le\frac34.
\]
It is not difficult to see that $\Wp^p(\mu_t,\nu)=2^p t$. Take $t_0=(\underline{\delta}/2)^p$ and $t_1=t_0+h$ with $h=(16\sqrt n)^{-1}$. Since $\underline{\delta}\le1$, we have $t_0\le1/2$ and $t_1\le9/16$. By Lemma \ref{lem: rootgap},
\[
\Wp(\mu_{t_1},\nu) - \Wp(\mu_{t_0},\nu) = (\underline{\delta}^p + 2^p h)^{1/p} -\underline{\delta} \asymp_p n^{-1/(2p)} \wedge (\underline{\delta}^{1-p}n^{-1/2}).
\]
By tensorization of $1+\chi^2$ and the Cauchy--Schwarz inequality,
\[
\begin{split}
 \TV(\mu_{t_1}^{\otimes n} \otimes \nu^{\otimes m},\mu_{t_0}^{\otimes n}\otimes \nu^{\otimes m}) &= \TV(\mu_{t_1}^{\otimes n},\mu_{t_0}^{\otimes n}) \le \frac{1}{2}\sqrt{\chi^2 (\mu_{t_1}^{\otimes n},\mu_{t_0}^{\otimes n})} \\
 &\le \frac{1}{2}\sqrt{(1+\chi^2(\mu_{t_1},\mu_{t_0}))^n-1} \\
 &\le \frac{1}{2}\sqrt{\left(1+\frac{h^2}{(1/4+t_0)(3/4-t_0)}\right)^n-1} \\
 &\le \frac{1}{2}\sqrt{e^{(16/3)nh^2}-1} = \frac{1}{2}\sqrt{e^{1/48}-1} < \frac{1}{2}.
 \end{split}
\]
Applying Le Cam's two-point method (cf. Chapter 2 in \cite{tsybakov2009introduction}), we obtain
\[
 \calR^{\Wp}_{n,m,1}(\underline{\delta})
 \gtrsim_p n^{-1/(2p)} \wedge (\underline{\delta}^{1-p}n^{-1/2}),
\]
as desired. 
\end{proof}

\begin{proposition}[Lower bound in $4<d\le2p$]
\label{prop: Wp appendix high lower}
Let $4 <d \le 2p$ and $n \le m$. For $0<\underline{\delta}\le 3/4$ and sufficiently large $n$ (depending only on $p,d$), we have
\[
 (\log n) \calR^{\Wp}_{n,m,d}(\underline{\delta})
 \gtrsim_{p,d} n^{-1/(2p)}\wedge
 \begin{cases}
 \big(\underline{\delta}^{1-p}n^{-1/2}+\underline{\delta}^{-1}n^{-2/d}\big), & p\le d, \\
 \big(\underline{\delta}^{1-p}n^{-1/2}+\underline{\delta}^{1-2p/d}n^{-2/d}\big),& p>d.
 \end{cases}
\]
\end{proposition}

\begin{proof}
If $\underline{\delta} \le n^{-1/(2p)}$, then the right-hand side of
(\ref{eq: Wp two-point}) reduces to $n^{-1/(2p)}$, and the claim follows from the previous proposition. Hence, we will assume that $\underline{\delta} > n^{-1/(2p)}$.

\medskip
\underline{Step 1} (Construction of random discrete measures).
We consider a similar construction to Section~\ref{sec: lower bound large p}. 
Let $\mathfrak{t}_z(y)=y+z$, and take
\[
 (z,\vartheta)=
 \begin{cases}
  (\underline{\delta}e_1,1),& p\le d,\\
  (3e_1/4,(4\underline{\delta}/3)^p),& p>d,
 \end{cases}
\]
which satisfies $\vartheta\|z\|^p=\underline{\delta}^p$. Recall that $\mathfrak{u}$ is the uniform distribution on $[M]$. Let $0<\eta\le1/10$ be chosen in Step~3 below. In view of Lemma \ref{lem: Wp testing}, we take
\begin{equation}
 k=n\wedge\lceil20n\vartheta\rceil \quad \text{and} \quad M=\lceil C_1\eta^{-1}k\log k\rceil.
 \label{eq: sample size choice}
\end{equation}
Let $F:[M] \to \calS$ be a uniform random bijection from Lemma \ref{lem: Wp packing} with $\calB=B_d(e_1/6,1/12)$.

Construct $(\mu_{F,\frakq},\overline{\nu})$ as in Section \ref{sec: lower bound large p},
\[
\mu_{F,\frakq} = \sum_{i=1}^M \frakq_i \frac{\delta_{F(i)}+\delta_{-F(i)}}{2}, \quad \overline{\nu} = \mu_{F,\fraku} = \frac{1}{2M} \sum_{x \in \calS \cup (-\calS)}\delta_{x},
\]
where $\frakq \in \calP([M])$. Now, we set
\[
 \breve{\mu}_{F,\frakq}=\vartheta \mu_{F,\frakq}+(1-\vartheta)\delta_{-e_1} \quad \text{and} \quad 
 \breve{\nu}=\vartheta(\mathfrak{t}_{z})_\#\overline{\nu}+(1-\vartheta)\delta_{-e_1},
 \]
 which are both supported in $B_d(0,1)$ for $\underline{\delta} \le 3/4$. 

We shall verify that 
\begin{equation}
\Wp^p(\breve{\mu}_{F,\frakq},\breve{\nu}) = \vartheta \Wp^p(\mu_{F,\frakq},(\mathfrak{t}_z)_\#\overline{\nu}),
\label{eq: translation}
\end{equation}
which is trivial when $\vartheta = 1$. For $\vartheta < 1$ (i.e., $z=3e_1/4$), apply Lemma \ref{lem: Wp block} with $\calX_0=\calS\cup(-\calS)$, $\calY_0=\calX_0+z$, and $\calX_1=\calY_1=\{-e_1\}$.  The distance between a point in $\calX_0$ and a point in $\calY_0$ is at most $5/4$, while the distances from $-e_1$ to $\calX_0$ and $\calY_0$ are at least $3/4$ and $3/2$, respectively. Hence,
\[
\overline m_{00}+\overline m_{11}\le(5/4)^p<(3/4)^p+(3/2)^p\le\underline m_{01}+\underline m_{10},
\]
which gives (\ref{eq: translation}).

Finally, recalling (\ref{eq: Wp translate-rem}), we obtain
\begin{equation}
 \breve{T}_{\frakq}(F):=\Wp^p(\breve{\mu}_{F,\frakq},\breve{\nu})
 =\underline{\delta}^p+\vartheta\inf_{\pi\in\Pi(\mu_{F,\frakq},\overline{\nu})}
       \int B_f(z+y-x,z)\,d\pi(x,y),
\label{eq: Wp bregman}
\end{equation}
where $f(x) = \| x \|^p$ and $B_f$ is its Bregman divergence (recall that $p > 2$). 
In particular, we have $\breve{T}_{\frakq}(F) \ge \underline{\delta}^p$. Set
\[
\Delta=\vartheta\|z\|^{p-2}M^{-2/d}.
\]

\medskip
\underline{Step 2} (Lower bound for $\breve{T}_{\frakq}(F)$ under $\Dplus$).
Pick any $\frakq \in \Dplus$. Arguing as in Step 2 in Section \ref{sec: lower bound large p}, we obtain 
\[
\breve{T}_{\frakq}(F) \ge \underline{\delta}^p + c_3 \Delta
\]
for some $c_3 \gtrsim_{p,d} 1$. 

\medskip
\underline{Step 3} (Upper bound for $\breve{T}_{\frakq}(F)$ under $\Dminus$).
Pick any $\frakq\in\Dminus$. To upper bound $\breve{T}_{\frakq}(F)$, we use Lemma \ref{lem: dyadic} with $(\mu,\nu) = (F_{\#}\frakq,F_{\#}\fraku)$ to construct a coupling $\pi$ (that depends on $F$). 
Given $F$, generate  $(U,V)$ from $\pi$. 
For an independent Rademacher variable $\varepsilon$, the conditional distribution of 
$(\varepsilon U,\varepsilon V)$ given $F$ gives a coupling for $(\mu_{F,\frakq},\overline{\nu})$ and $\| \varepsilon U - \varepsilon V\| = \|U-V\|$. 
As such, combining (\ref{eq: Wp bregman}) and (\ref{eq: bregman}), we obtain
\begin{equation}
 0\le \E_F[\breve{T}_{\frakq}(F)-\underline{\delta}^p]\lesssim_p
 \vartheta\big(\|z\|^{p-2}\E[\|U-V\|^2]+\E[\|U-V\|^p]\big).
 \label{eq: trial cost}
\end{equation}

Consider the dyadic partition $\calQ_\ell$. For any $Q \in \calQ_{\ell}$, we observe (cf. Lemma 7 in \cite{nilesweed2022estimation}),
\[
\begin{split}
&\E_F[(F_{\#}\frakq)(Q)] = (F_{\#}\fraku)(Q), \\ 
&\Var_F((F_{\#}\frakq)(Q)) = \frac{(F_{\#}\fraku)(Q)(1-(F_{\#}\fraku)(Q))\chi^2(\frakq,\fraku)}{M-1}, \\ 
&\sum_{Q \in \calQ_\ell} |(F_{\#}\frakq)(Q) - (F_{\#}\fraku)(Q)| \le 2\TV(\frakq,\fraku). 
\end{split}
\]
The final inequality follows by  Scheff\'{e}'s inequality (cf. (\ref{eq: scheffe})). 
Recalling that $\chi^2(\frakq,\fraku) \le 9$ and $\TV(\frakq,\fraku) \le \eta$, and applying the Cauchy--Schwarz inequality, we obtain
\[
 \E_F \left [ \sum_{Q \in \calQ_\ell} |(F_{\#}\frakq)(Q) - (F_{\#}\fraku)(Q)| \right]
 \lesssim \eta \wedge \sum_{Q\in\calQ_\ell}\sqrt{\frac{(F_{\#}\fraku)(Q)}{M}}
 \lesssim \eta\wedge\sqrt{\frac{2^{d\ell}}{M}}.
\]
Now, applying Lemma \ref{lem: dyadic} and using (\ref{eq: dyadic series}), we have
\[
 \E[\|U-V\|^2] \lesssim_d\eta^{1-4/d} M^{-2/d},\quad
 \E[\|U-V\|^p] \lesssim_{p,d} M^{-1/2} (\log M)^{\ind_{\{ d=2p \}}}.
\]
For $p>d$, $\|U-V\|\le1/6$ and $\|z\|=3/4$, so that
$\|U-V\|^p\lesssim_p\|z\|^{p-2}\|U-V\|^2 \lesssim_p \|U-V\|^2$. For $p\le d$, we use $\vartheta=1$.
Substituting these bounds into (\ref{eq: trial cost}) and recalling
$\Delta=\vartheta\|z\|^{p-2}M^{-2/d}$, we obtain
\begin{equation}
 \E_F[\breve{T}_{\frakq}(F)-\underline{\delta}^p]\le
 C_3\left(\eta^{1-4/d}\Delta+\ind_{\{p\le d\}} M^{-1/2}(\log M)^{\ind_{\{d=2p\}}}\right),
 \label{eq: Wp upper}
\end{equation}
where $C_3 \lesssim_{p,d} 1$. 

Choose $0 < \eta \le 1/10$ (that may depend on $(p,d)$) so that
\[
C_3\eta^{1-4/d}\le \frac{c_3}{40},
\]
where $c_3$ is given in Step 2. 
For $p\le d$, $M\asymp_{p,d} n\log n$ and
\[
 M^{-1/2}(\log M)^{\ind_{\{d=2p\}}} \lesssim_{p,d}  n^{-1/2}\sqrt{\log n}.
\]
If in addition $n^{-1/2}\sqrt{\log n} \gtrsim_{p,d} \Delta$, then
Proposition \ref{prop: Wp two point} and $\underline{\delta}>n^{-1/(2p)}$ imply
\begin{equation}
 \underline{\delta}^{1-p}\Delta
 \lesssim_{p,d}\sqrt{\log n}\,\underline{\delta}^{1-p}n^{-1/2}
 \lesssim_{p,d}\sqrt{\log n}\,\calR^{\Wp}_{n,m,d}(\underline{\delta}).
 \label{eq: absorbed}
\end{equation}

\medskip
\underline{Step 4} (Application of Proposition~\ref{prop: fuzzy}).
Suppose $p>d$ or $n^{-1/2}\sqrt{\log n} \lesssim_{p,d} \Delta$. By adjusting $\eta$ if necessary, 
by  (\ref{eq: Wp upper}),
we obtain 
\[
\sup_{\frakq\in\Dminus}\E_F[\breve{T}_{\frakq}(F)-\underline{\delta}^p]\le \frac{c_3\Delta}{20}.
\]
Let $\Lambda_{\frakq}$ be the distribution of $(\breve{\mu}_{F,\frakq},\breve{\nu})$ on $\calP(\calX) \times \calP(\calX)$. Set
\[
 c=(\underline{\delta}^p+c_3\Delta/2)^{1/p} \quad \text{and} \quad
s = \frac{(\underline{\delta}^p+c_3\Delta)^{1/p}-(\underline{\delta}^p+c_3\Delta/2)^{1/p}}{2}.
\]
By Markov's inequality, for any $q \in \Dminus$, 
\[
\Lambda_{\frakq}(\Wp>c) \le\frac{c_3\Delta/20}{c_3\Delta/2}=0.1.
\]
On the other hand, from Step 2, for any $\frakq \in \Dplus$,
\[
 \Lambda_{\frakq}(\Wp \ge c+2s)=1.
\]

Pick any estimator $\hat W$ and put $\psi=\ind_{\{\hat W\ge c+s\}}$.
Generate $I_1,I_2,\dots \stackrel{\text{i.i.d.}}{\sim} \frakq$, $J_1,\dots,J_n \stackrel{\text{i.i.d.}}{\sim} \mathsf{Bern}(\vartheta)$, $\varepsilon_1,\dots,\varepsilon_n \stackrel{\text{i.i.d.}}{\sim} \mathsf{Rad}$, independently of one another\footnote{We use $\mathsf{Bern}$ and $\mathsf{Rad}$ to denote the Bernoulli and Rademacher distributions, respectively.}, and independently from $F$. 
Generate
\[
 X_i=\begin{cases}
 -e_1, &J_i=0, \\
 \varepsilon_i F(I_{J_1+\cdots+J_i}), & J_i=1,
 \end{cases}
\]
and $Y_1,\dots,Y_m \stackrel{\text{i.i.d.}}{\sim} \breve{\nu}$ independently from everything else. Given $F$, 
\[
(X_1,\dots,X_n,Y_1,\dots,Y_m) \sim \breve{\mu}_{F,\frakq}^{\otimes n} \otimes \breve{\nu}^{\otimes m}.
\]
Recall that we have chosen $k$ and $M$ as in (\ref{eq: sample size choice}). In this case,  $k$ may be strictly smaller than $n$, and the test $\psi$ depends on extra $I_{k+1},\dots,I_{n}$ that depend on the parameter $\frakq$. Hence, one may not directly apply Lemma \ref{lem: Wp testing} to $\psi$. As such, we consider a modified test
\[
\psi_k = \psi \ind_{\{ \sum_{i=1}^n J_i \le k\}},
\]
which only uses $I_1,\dots,I_k$. For any $q\in\calP([M])$, 
\[
\Prob_{\frakq}(\psi \ne \psi_k) \le \Prob \left ( \sum_{i=1}^n J_i > k\right ) \le \frac{n\vartheta}{k} \le 0.05.
\]
Since $k \asymp n\vartheta \gtrsim_p \sqrt{n}$ (recall that $\underline{\delta} > n^{-1/(2p)}$), by Lemma \ref{lem: Wp testing}, for sufficiently large $n$, one can find 
$\frakq^-\in\Dminus$ and $\frakq^+\in\Dplus$ such that
\[
 \Prob_{\frakq^-}(\psi_k=1)+\Prob_{\frakq^+}(\psi_k=0)
 \ge 0.8,
\]
which implies
\[
 \Prob_{\frakq^-}(\psi=1)+\Prob_{\frakq^+}(\psi=0)
 \ge 0.7.
\]

We shall apply Proposition \ref{prop: fuzzy} with $\sfX=\calX^{n+m}, \Theta = \{ (\mu,\nu) \in \calP(\calX) \times \calP(\calX): \Wp(\mu,\nu) \ge \underline{\delta} \}, P_\theta = \mu^{\otimes n} \otimes \nu^{\otimes m}$,  $T(\theta) = \Wp(\mu,\nu)$ for $\theta = (\mu,\nu)$, $\Lambda_0 = \Lambda_{\frakq^-}$, and $\Lambda_1 = \Lambda_{\frakq^+}$. 
We have already verified Conditions (i) and (ii) with $\beta_0=0.1, \beta_1=0$, and $\alpha=0.3$. We conclude that
\[
\sup_{\Wp(\mu,\nu) \ge \underline{\delta}}\Prob_{\mu,\nu}(|\hat W-\Wp(\mu,\nu)|\ge s) \ge 0.3. 
\]
We have established $\calR^{\Wp}_{n,m,d}(\underline{\delta}) \gtrsim_{p,d} s \asymp_{p,d} \Delta^{1/p} \wedge (\underline{\delta}^{1-p}\Delta)$ when $p>d$ or $n^{-1/2}\sqrt{\log n} \lesssim_{p,d} \Delta$. 
We shall verify $\Delta^{1/p} \gtrsim_{p,d} \underline{\delta}^{1-p}\Delta$, which implies $s \asymp_{p,d} \underline{\delta}^{1-p}\Delta$. Observe that $\Delta^{1/p} \gtrsim_{p,d} \underline{\delta}^{1-p}\Delta$ is equivalent to $\Delta \lesssim_{p,d} \underline{\delta}^p$. Recall $\underline{\delta} > n^{-1/(2p)}$. For $p \le d$,
\[
\Delta/\underline{\delta}^p = \underline{\delta}^{-2} M^{-2/d} \lesssim_{p,d} n^{1/p} (n \log n)^{-2/d} \lesssim 1.
\]
For $p > d$, $\Delta/\underline{\delta}^p \lesssim M^{-2/d} \lesssim 1$. We conclude that $\calR^{\Wp}_{n,m,d}(\underline{\delta}) \gtrsim_{p,d} \underline{\delta}^{1-p}\Delta$ when $p>d$ or $n^{-1/2}\sqrt{\log n} \lesssim_{p,d} \Delta$. 

When $p \le d$ and $n^{-1/2}\sqrt{\log n} \gtrsim_{p,d} \Delta$, we have already established (\ref{eq: absorbed}). 
Finally, recalling that $\vartheta = (\underline{\delta}/\|z\|)^p$, we have 
\[
 \begin{split}
 \underline{\delta}^{1-p}\Delta
 &\asymp_{p,d}\underline{\delta}^{1-p}\vartheta^{1-2/d}\|z\|^{p-2} (n\log n)^{-2/d} \\
 &=\underline{\delta}^{1-2p/d}\|z\|^{2p/d-2}(n\log n)^{-2/d}\\
 &\asymp_{p,d}(n\log n)^{-2/d}
 \begin{cases}
 \underline{\delta}^{-1}, & p\le d,
 \\ \underline{\delta}^{1-2p/d}, & p>d.
 \end{cases}
 \end{split}
\]
Combining the preceding proposition\footnote{Use the simple fact that $R \ge a \vee b \asymp a+b$ if $R \ge a$ and $R \ge b$ for some $a,b > 0$.}, we obtain the result. 
\end{proof}

\subsection{Upper bound}

By symmetry, we assume $N=n \le m$. Pick any $(\mu,\nu)$ with $t := \Wp(\mu,\nu) \ge \underline{\delta}$ and let $\hat{W} = \Wp(\hat{\mu}_n,\hat{\nu}_m)$. 
We first observe that, by Theorem 1 in \cite{fournier2015on}, 
\begin{equation}
\E[|\hat{W}-t|] \lesssim_{p,d} n^{-1/(2p)} (\log n)^{\ind_{\{d=2p\}}/p}.
\label{eq: first bound}
\end{equation}

We use the notation of Section \ref{sec: power bound}. In Steps 1--3 of that section, the $d > 2p$ condition was used only to invoke $\E[S]\lesssim_{p,d}n^{-p/d}$ and $n^{-1/2}\le n^{-p/d}$. With that in mind, applying \eqref{eq: Wp cost-bias} and \eqref{eq: bias} and taking into account
Remark~\ref{rem: upper bound for small p} and Lemma~\ref{lem: bregman}, we obtain
\begin{equation}
 \E[|\hat W^p-t^p|]\lesssim_p n^{-1/2}+
 \begin{cases}
\E\left[n^{-1}\sum_{i=1}^n\|b_i\|^p\right], & 1\le p\le2,\\
\E\left[n^{-1}\sum_{i=1}^n
       \big(\|a_i\|^{p-2}\|b_i\|^2+\|b_i\|^p\big)\right], &p>2.
 \end{cases}
 \label{eq: Wp permutation bound}
\end{equation}
Observe that $\sigma$ was chosen as a minimizer of (\ref{eq: permutation}), but can be replaced with any permutation of $[n]$ that may depend on $(Z_1,\dots,Z_n)$ and $(Z_1',\dots,Z_n')$. We will consider the following two cases: 1) $d \le 4$ and 2) $4 <d \le 2p$.  

\medskip

\underline{Step 1} (Case 1). 
For any $r \ge 1$, choose $\sigma$ that minimizes $n^{-1}\sum_{i=1}^n\|b_i\|^r = n^{-1}\sum_{i=1}^n \|Z_i-Z_{\sigma(i)}'\|^r$. By Theorem 1 in \cite{fournier2015on},
\[
\E\left [\frac{1}{n}\sum_{i=1}^n\|b_i\|^r\right] \lesssim_{r,d} n^{-\frac{r}{d \vee (2r)}} (\log n)^{\ind_{\{ d=2r \}}}.
\]
Using $\|a_i\|^{p-2}\|b_i\|^2+\|b_i\|^p \lesssim_p \| b_i \|^2$ for $p > 2$, we obtain, from (\ref{eq: Wp permutation bound}), 
\[
\E[|\hat W^p-t^p|] \lesssim_{p,d} n^{-\frac{\alpha}{d \vee (2\alpha)}} (\log n)^{\ind_{\{ d=2\alpha \}}},
\]
where $\alpha = p \wedge 2$, 
which implies
\begin{equation}
\E[|\hat W-t|] \lesssim_{p,d} \underline{\delta}^{1-p} n^{-\frac{\alpha}{d \vee (2\alpha)}} (\log n)^{\ind_{\{ d=2\alpha \}}}.
\label{eq: second bound}
\end{equation}
Taking the minimum of (\ref{eq: first bound}) and (\ref{eq: second bound}) covers $d \le 4$.
\medskip

\underline{Step 2} (Case 2). Suppose $4<d\le2p$, which entails $p>2$.
Let $D:\R^d\to\R^d$ be a Lipschitz function constructed from Lemma \ref{lem: shortening}. Recall that $(U_i,V_i)\sim\pi$ and $Z_i=(U_i+V_i)/2$. For $z,z'\in\mathsf{D}$, let
\[
 \sfc(z,z')=(\|D(z)\|\vee\|D(z')\|)^{p-2}\|z-z'\|^2+\|z-z'\|^p.
\]
Let $\hat{\gamma}_n,\hat{\gamma}_n'$ be the empirical distributions of
$\{Z_i\}_{i=1}^n,\{Z_i'\}_{i=1}^n$, respectively. Observe that $D(Z_i) = U_i-V_i=\Delta_i$. 
From (\ref{eq: Wp permutation bound}) and the matching representation of the OT cost for empirical distributions (cf. (\ref{eq: matching})), we have
\[
 \E[|\hat W^p-t^p|]\lesssim_p n^{-1/2}+\E[T_{\sfc}(\hat{\gamma}_n,\hat{\gamma}_n')].
\]

We use Lemma \ref{lem: dyadic} to bound $\E[T_{\sfc}(\hat{\gamma}_n,\hat{\gamma}_n')]$.
For $Q\in\calQ_\ell$, let $r_Q=\Prob(Z_1\in Q), w_Q=\inf_{z\in Q}\|D(z)\|^{p-2}$,
and $a_Q=|\hat{\gamma}_n(Q)-\hat{\gamma}_n'(Q)|$. Since $D$ is $2\sqrt{p-1}$-Lipschitz and
$\operatorname{diam}(Q)\lesssim_d 2^{-\ell}$,
\[
 \overline{\sfc}_Q
 \lesssim_{p,d} 2^{-2\ell}\big(w_Q^{1/(p-2)}+2^{-\ell}\big)^{p-2}+2^{-p\ell}
 \lesssim_{p,d} 2^{-2\ell}w_Q+2^{-p\ell}.
\]
Observe that
$\E[a_{Q'}]\le 2r_{Q'}\wedge\sqrt{2r_{Q'}/n}$ and $\sum_{Q'\prec Q}\E[a_{Q'}]\lesssim_d r_Q \wedge \sqrt{r_Q/n}$.
Applying \eqref{eq: dyadic moment bound}, we obtain
\[
 \E[T_{\sfc}(\hat{\gamma}_n,\hat{\gamma}_n')]
 \lesssim_{p,d}\sum_{\ell=0}^\infty\sum_{Q\in\calQ_\ell}
       (2^{-2\ell}w_Q+2^{-p\ell})\big( r_Q \wedge \sqrt{r_Q/n} \big ).
\]
By (\ref{eq: dyadic series}) and $d\le2p$,
\[
 \sum_{\ell=0}^\infty 2^{-p\ell}\sum_{Q\in\calQ_\ell}\big ( r_Q \wedge \sqrt{r_Q/n} \big)
 \le\sum_{\ell=0}^\infty 2^{-p\ell}\big(1\wedge\sqrt{2^{d\ell}/n}\big)
 \lesssim_{p,d}n^{-1/2}(\log n)^{\ind_{\{d=2p\}}}.
\]

Put $q=d/(d-2)$ and $B=(\E[\|\Delta_1\|^{(p-2)q}])^{1/q}$.
Observe that $\sum_{Q\in\calQ_\ell}w_Q^q r_Q\le B^q$. Since $d>4$,
$r \wedge \sqrt{r/n} \le r\wedge(n^{-2/d}r^{1/q})$ for $0 < r \le 1$. By H\"older's inequality,
\[
 \sum_{Q\in\calQ_\ell}w_Q \big ( r_Q \wedge \sqrt{r_Q/n} \big)
 \le B\big(1\wedge n^{-2/d}|\calQ_\ell|^{1-1/q}\big)
 =B\big(1\wedge n^{-2/d}2^{2\ell}\big).
\]
Splitting  at $\lfloor(\log_2 n)/d\rfloor$, we obtain
\[
 \sum_{\ell=0}^\infty 2^{-2\ell}\sum_{Q\in\calQ_\ell}w_Q\big ( r_Q \wedge \sqrt{r_Q/n} \big)
 \le B\sum_{\ell=0}^\infty(2^{-2\ell}\wedge n^{-2/d})
 \lesssim_d Bn^{-2/d}\log n .
\]

Finally, $\E[\|\Delta_1\|^p]=t^p$, $\|\Delta_1\|\le2$, and $(p-2)q\le p$
if and only if $p\le d$. Hence, when $p \le d$,
\[
B^q = \E[\|\Delta_1\|^{(p-2)q}] \le (\E[\|\Delta_1\|^p])^{(p-2)q/p} \le t^{(p-2)q},
\]
while when $p > d$,
\[
B^q \le 2^{(p-2)q-p} \E[\|\Delta_1\|^p] \lesssim t^{p}.
\]
Summarizing, we obtain
\[
 B\lesssim_{p,d}
 \begin{cases}
 t^{p-2},&p\le d,\\
 t^{p/q},&p>d,
 \end{cases}
\]
which implies 
\[
\E[|\hat{W}^p-t^p|] \lesssim_{p,d} n^{-1/2}(\log n)^{\ind_{\{ d=2p \}}} + 
\begin{cases}
t^{p-2} n^{-2/d}\log n, & p \le d, \\
t^{p/q} n^{-2/d}\log n, & p > d.
\end{cases}
\]
Using $|\hat W-t|\le t^{1-p}|\hat W^p-t^p|$ and $t\ge\underline{\delta}$, we conclude that
\[
 \E[|\hat W-t|]\lesssim_{p,d} \underline{\delta}^{1-p}n^{-1/2}(\log n)^{\ind_{\{d=2p\}}} + 
 \begin{cases}
\underline{\delta}^{-1}n^{-2/d}\log n,& p\le d,\\
\underline{\delta}^{1-2p/d}n^{-2/d}\log n,& p>d.
 \end{cases}
\]
Taking the minimum with \eqref{eq: first bound} completes the proof.
\qed

\section{Proof of Theorem \ref{thm: GW main}}
\label{sec: GW proof}

\subsection{Preliminaries}
The lower-bound proof adapts several elements from the proof of Theorem 3.2 in \cite{kato2025convergence} and the lower-bound proof for $\Wp$ with $1 \le p \le 2$ in Theorem \ref{thm: Wp main}. One key difference from the $\Wp$ case, because of the lack of a decomposition formula analogous to Lemma \ref{lem: Wp block},  is that instead of adding atoms to $(\mu_{F,\frakq},\overline{\nu})$, we dilate $\overline{\nu}$ and account for the change in $\GW$. The next lemma concerns the effect of dilation on $\GW$. 

\begin{lemma}
\label{lem: GW-scaling}
Let $\mu \in \calP(\R^{d_x})$ and $\nu \in \calP(\R^{d_y})$ be compactly supported. For $k\in\{d_x,d_y\}$ and $c\ge0$, let $\mathfrak{s}_c:\R^k\to\R^k$ denote the dilation $\mathfrak{s}_c(z)=cz$. Then for any $a,b \ge 0$,
\[
\GW^2((\mathfrak{s}_a)_{\#}\mu,(\mathfrak{s}_b)_{\#}\nu)
 =a^2b^2 \GW^2(\mu,\nu) +(a^2-b^2)\big( a^2\sfA(\mu)-b^2\sfA(\nu) \big).
\]
\end{lemma}

\begin{proof}
It suffices to consider the case where $ab > 0$. 
Observe that the pushforward under $(x,y) \mapsto (\mathfrak{s}_a(x),\mathfrak{s}_b(y))$ is a bijection from $\Pi(\mu,\nu)$ onto $\Pi((\mathfrak{s}_a)_{\#}\mu,(\mathfrak{s}_b)_{\#}\nu)$. Using the expression (\ref{eq: GW-fourth-moments}), we have 
\[
\begin{split}
\GW^2((\mathfrak{s}_a)_{\#}\mu,(\mathfrak{s}_b)_{\#}\nu) &= a^4\sfA(\mu) + b^4\sfA(\nu) \\
&\quad -2a^2b^2 \sup_{\pi \in \Pi(\mu,\nu)} \iint \|x-x'\|^2\|y-y'\|^2 \, d\pi(x,y)d\pi(x',y'). 
\end{split}
\]
Subtracting $a^2b^2\GW^2(\mu,\nu)$ from $\GW^2((\mathfrak{s}_a)_{\#}\mu,(\mathfrak{s}_b)_{\#}\nu)$ yields
\[
a^4\sfA(\mu) + b^4\sfA(\nu) - a^2b^2 \big ( \sfA(\mu) + \sfA(\nu)\big ) = (a^2-b^2)\big( a^2\sfA(\mu)-b^2\sfA(\nu) \big),
\]
which gives the result. 
\end{proof}

We will apply Lemma \ref{lem: GW-scaling} to $(\mu,\nu) = (\mu_{F,\frakq},\overline{\nu})$ and $(a,b)=(1,b)$ where $0<b<1$ is chosen so that $(1-b^2)^2\sfA(\overline{\nu})=\underline{\delta}^2+\Delta$, with $\Delta>0$ specified below. 
Still, one has to account for the fluctuation of $\sfA(\mu_{F,q})-\sfA(\overline{\nu})$, which is the purpose of the next lemma. 

\begin{lemma}
\label{lem: A-stability}
Let $\mathcal{S} = \{ x_1,\dots,x_M \} \subset B_d(0,1)$ be any set of $M\ge2$ points. Let $\mathfrak{u}$ be the uniform distribution on $[M]$ and let $q \in \calP([M])$ be such that $\chi^2(\mathfrak{q},\mathfrak{u})\le9$. If $F$ is a uniform random bijection from $[M]$ onto $\mathcal{S}$, then $\mu_{F,\frakq}=F_{\#}\frakq$ and $\overline{\nu}=F_{\#}\fraku$ satisfy
\begin{equation}\label{eq:A-stability}
 \E_{F}[|\sfA(\mu_{F,q})-\sfA(\overline{\nu})|^2]\lesssim 1/M.
\end{equation}
\end{lemma}

\begin{proof}
Let $V=V(F) = \sfA(\mu_{F,\frakq}), \xi=\chi^2(q,u)$, and $\frakq_i=\frakq(\{i\})$. It suffices to consider the case when $\xi > 0$. To bound the second moment, we shall bound the variance and the squared bias. Observe that $(F(1),\dots,F(M))$ follows sampling without replacement from $\{ x_1,\dots,x_M \}$.

To bound $\Var_F(V)$, let $H_{ij}=\|F(i)-F(j)\|^4$ for $(i,j)\in[M]^2$. We have 
\[
V (F)= \sum_{i,j=1}^M\frakq_i \frakq_j H_{ij}  = \sum_{i \ne j} \frakq_i \frakq_j H_{ij}. 
\]
If $\tau_{ij}$ is a transposition between $i$ and $j$, then for $i \ne j$,
\[
 V(F)-V(F\circ\tau_{ij})
 =2(\frakq_i-\frakq_j)\sum_{k\in[M]\setminus\{i,j\}} \frakq_k(H_{ik}-H_{jk}).
\]
Since $0\le H_{ij}\le16$, the absolute difference is at most $32|\frakq_i-\frakq_j|$.  Apply Proposition 1 in \cite{sambale2022concentration} to $V$ and $-V$. Since 
\[
\sum_{i < j}(\frakq_i-\frakq_j)^2 = \frac{1}{2}\sum_{i,j}(\frakq_i-\frakq_j)^2 = M \sum_{i} \frakq_i^2 -1 = \chi^2(\frakq,\fraku) = \xi,
\]
we obtain 
\[
\Prob_F(|V-\E_{F}[V]|\ge t)\le 2e^{-cMt^2/\xi},
\]
for $t>0$, where $c \gtrsim 1$. Integrating gives 
\[
\Var_F (V) = \int_0^\infty 2t\Prob_F(|V-\E_{F}[V]|\ge t) \, dt  \lesssim \xi/M.
\]

For $i\ne j$, the pair $(F(i),F(j))$ is uniformly distributed over the $M(M-1)$ pairs of distinct points from $\calS$, so 
\[
\E_F[H_{ij}] = \frac{1}{M(M-1)} \sum_{i \ne j}\|x_i-x_j\|^4 = \frac{M}{M-1} \sfA(\overline{\nu}),
\]
where we used $\sfA(\overline{\nu}) = M^{-2}\sum_{i \ne j} \|x_i-x_j\|^4$. Hence,
\[
\begin{split}
 \E_{F}[V]&=\left ( \sum_{i \ne j} \frakq_i \frakq_j \right) \frac{M\sfA(\overline{\nu})}{M-1}  =\left(1-\sum_{i=1}^M \frakq_i^2\right)\frac{M \sfA(\overline{\nu})}{M-1} \\
 &=\left(1-\frac{\xi}{M-1}\right)\sfA(\overline{\nu}),
 \end{split}
\]
where we used $\sum_{i=1}^M \frakq_i^2 = (1+\xi)/M$. 
Adding the squared bias 
\[
\left (\frac{\xi \sfA(\overline{\nu})}{M-1} \right )^2 \lesssim \frac{1}{M}
\]
to the variance bound proves~(\ref{eq:A-stability}), where we used $\xi \le 9$ and $\sfA(\overline{\nu}) \le 16$.
\end{proof}

\subsection{Lower bound}

Recall $d = d_x \wedge d_y$.  By Lemma \ref{lem: minimax risk}~(ii), it suffices to consider the case where $d_x=d_y = d$. In addition, by symmetry, $\calR^{\mathsf{GW}}_{n,m,d,d} = \calR^{\mathsf{GW}}_{m,n,d,d}$, so that we may further assume that $N=n \le m$. The lower bound follows by combining Propositions \ref{prop: GW prelim lower} and \ref{prop: GW lower large d} below. 

\begin{proposition}[Preliminary lower bound]
\label{prop: GW prelim lower}
Under $n \le m$, for any $0 < \underline{\delta} \le 1$, we have
\[
\calR^{\mathsf{GW}}_{n,m,d,d}(\underline{\delta}) \gtrsim n^{-1/4} \wedge (\underline{\delta}^{-1}n^{-1/2}). 
\]
\end{proposition}
\begin{proof}
The proof is an adaptation of Step 1 in the proof of Theorem 3.2 in \cite{kato2025convergence}.
By Lemma \ref{lem: minimax risk}~(iii), it suffices to consider the $d=1$ case.
Let $\nu$ be the uniform distribution on $\{-1,1\}$, and put
\[
 \mu_t=\left ( \frac{1}{2} + t \right )\delta_{-1}+\left ( \frac{1}{2} - t \right )\delta_1,\quad 0\le t\le \frac{1}{2}.
\]
As in Lemma 4.7 in \cite{kato2025convergence}, we have
\[
 \GW^2(\mu_t,\nu)=32t(1-t).
\]

Since $\GW^2(\mu_{1/16},\nu)=15/8>1$, we may choose $0<t_0<1/16$ so that $\GW(\mu_{t_0},\nu)=\underline{\delta}$. Put $h=(16\sqrt n)^{-1}$ and $t_1=t_0+h<1/8$. By construction,  $\GW^2(\mu_{t_1},\nu) - \GW^2(\mu_{t_0},\nu) = 32h(1-t_0-t_1)\asymp h$. By Lemma \ref{lem: rootgap},
\[
 \GW(\mu_{t_1},\nu)-\GW(\mu_{t_0},\nu)\asymp\sqrt h\wedge(\underline{\delta}^{-1}h) \asymp n^{-1/4}\wedge(\underline{\delta}^{-1}n^{-1/2}).
\]
 By tensorization of $1+\chi^2$ and the Cauchy--Schwarz inequality,
\[
\begin{split}
 \TV(\mu_{t_1}^{\otimes n} \otimes \nu^{\otimes m},\mu_{t_0}^{\otimes n}\otimes \nu^{\otimes m}) &= \TV(\mu_{t_1}^{\otimes n},\mu_{t_0}^{\otimes n}) \le \frac{1}{2}\sqrt{\chi^2 (\mu_{t_1}^{\otimes n},\mu_{t_0}^{\otimes n})} \\
 &\le \frac{1}{2}\sqrt{(1+\chi^2(\mu_{t_1},\mu_{t_0}))^n-1} \\
 &\le \frac{1}{2}\sqrt{\left(1+\frac{h^2}{1/4-t_0^2}\right)^n-1} \\
 &\le  \frac{1}{2}\sqrt{e^{5nh^2}-1} =  \frac{1}{2}\sqrt{e^{5/16^2}-1} < \frac{1}{2}.
 \end{split}
\]
Applying Le Cam's two-point method (cf. Chapter 2 in \cite{tsybakov2009introduction}), we obtain
\[
 \calR^{\mathsf{GW}}_{n,m,1,1}(\underline{\delta})
 \gtrsim n^{-1/4}\wedge(\underline{\delta}^{-1}n^{-1/2}),
\]
as desired. 
\end{proof}
\begin{proposition}[Lower bound in $d > 4$]
\label{prop: GW lower large d}
Let $d > 4, 0 < \underline{\delta} \le 1$, and $n \le m$. For sufficiently large $n$ (depending only on $d$), we have
\[
\calR^{\mathsf{GW}}_{n,m,d,d}(\underline{\delta}) \gtrsim_d (n\log n)^{-1/d}\wedge\big(\underline{\delta}^{-1}(n\log n)^{-2/d}\big). 
\]
\end{proposition}

In the proof below, for any $\rho \in \calP(\R^d)$ with finite second moment, let $\Cov (\rho)$ denote the covariance matrix of $\rho$. We use $\lambda_{\min} (B)$ to denote the smallest eigenvalue of a symmetric matrix $B$. 

\begin{proof}
For ease of notation, set $\calX = B_d(0,1)$. We shall apply Proposition \ref{prop: fuzzy} with $\sfX=\calX^{n+m}, \Theta = \{ (\mu,\nu) \in \calP(\calX) \times \calP(\calX): \GW(\mu,\nu) \ge \underline{\delta} \}, P_\theta = \mu^{\otimes n} \otimes \nu^{\otimes m}$, and $T(\theta) = \GW(\mu,\nu)$ for $\theta = (\mu,\nu)$.

\medskip

\underline{Step 1} (Construction of random discrete measures).
Let $M$ be a sufficiently large positive integer to be determined later. 
 As in the lower-bound proof for $1 \le p \le 2$ in Theorem \ref{thm: Wp main}, we consider, for any $q \in \calP([M])$,
 \[
 \mu_{F,\frakq}=F_{\#}\frakq \quad \text{and} \quad \overline{\nu} = F_{\#}\fraku,
 \]
and suitably dilate $\overline{\nu}$. Here, $F$ is a uniform random bijection from $[M]$ onto a set $\calS\subset B_d(0,1)$. We choose $\calS$ so that the minimum spacing in $\calS$ is $\gtrsim_d M^{-1/d}$, $\overline{\nu}$ has mean zero, $\lambda_{\min}(\Cov(\overline{\nu}))\gtrsim_d 1$, and  $\sfA(\overline{\nu}) \ge 2$. 

To see that such $\calS$ exists,
let $k=\lfloor M/2\rfloor$. By translating and rescaling the points in Lemma 4.6 of \cite{kato2025convergence}, we may choose $x_1,\ldots,x_k$ in $B_d(7e_1/8,1/8)$ at pairwise distances $\gtrsim_d k^{-1/d}$ such that the smallest eigenvalue of the covariance matrix of  $\rho:=k^{-1}\sum_{i=1}^k \delta_{x_i}$ is $\gtrsim_d1$. Take $\calS=\{\pm x_1,\ldots,\pm x_k\}$ and include the origin if $M$ is odd. Since $\|x_i+x_j\|\ge3/2$ for all $i,j$ and $\|x_i\|\ge3/4$, the minimum spacing in $\calS$ is $\gtrsim_d M^{-1/d}$. By symmetry, $\overline{\nu} := |\calS|^{-1}\sum_{x \in \calS}\delta_{x}$ has mean zero, and  $ \lambda_{\min}(\Cov (\overline{\nu})) \ge (2k/M) \lambda_{\min}(\Cov (\rho)) \gtrsim_d 1$. Also, taking one point from each of the two reflected sets in the definition of $\sfA$, we see that
\[
 \sfA(\overline{\nu}) \ge 2\left(\frac{k}{M}\right)^2\left(\frac{3}{2}\right)^4 \ge 2
\]
for sufficiently large $M$.

Recall the definitions of $\Dplus$ and $\Dminus$ before Lemma \ref{lem: Wp testing}.
The preceding construction of $\calS$ allows us to apply the comparison in equation (29) of \cite{kato2025convergence}. As such,  there exist positive constants $c_4, C_4$ that depend only on $d$ such that
\[
 \begin{split}
 &\Prob_F \big ( \GW(\mu_{F,\frakq},\overline{\nu}) \le C_4 M^{-1/d}\eta^{1/2-2/d} \big) \ge 0.9, \quad \frakq \in \Dminus, \\
 &\Prob_F \big ( \GW(\mu_{F,\frakq},\overline{\nu}) \ge c_4 M^{-1/d} \big) = 1, \quad \frakq \in \Dplus.
 \end{split}
\]
Since $d>4$, one may choose sufficiently small $0<\eta\le1/10$ (that may depend on $d$) such that 
\[
C_4\eta^{1/2-2/d}\le \frac{c_4}{4}.
\]
Take $M$ as in (\ref{eq: M-choice}) with $k=n$ and set 
\[
\Delta=\frac{c_4^2M^{-2/d}}{16} \asymp_d (n \log n)^{-2/d}. 
\]
The choice of $\Delta$ guarantees that
\begin{equation}
 \begin{split}
 &\Prob_F \big ( \GW^2(\mu_{F,\frakq},\overline{\nu}) \le \Delta \big) \ge 0.9, \quad \frakq \in \Dminus, \\
 &\Prob_F \big ( \GW^2(\mu_{F,\frakq},\overline{\nu}) \ge 16 \Delta \big) = 1, \quad \frakq \in \Dplus.
 \end{split}
 \label{eq: step 1}
\end{equation}

\medskip

\underline{Step 2} (Separation after dilation).
Choose $0<b<1$ so that $(1-b^2)^2\sfA(\overline{\nu})=\underline{\delta}^2+\Delta$. Since $\sfA(\overline{\nu})\ge2$ and $\underline{\delta}\le1$, we have $b^2\ge1/4$ once $\Delta\le1/8$. 
By Lemma \ref{lem: GW-scaling}, we have
\[
 \begin{split}
 \GW^2(\mu_{F,\frakq},(\mathfrak{s}_b)_{\#}\overline{\nu})
 &=\underline{\delta}^2+\Delta+b^2\GW^2(\mu_{F,\frakq},\overline{\nu})\\
 &\quad +(1-b^2)\big(\sfA(\mu_{F,\frakq})-\sfA(\overline{\nu})\big).
 \end{split}
\]
By Lemma \ref{lem: A-stability} and Markov's inequality,
\[
 \sup_{\frakq: \chi^2(\frakq,\fraku) \le 9} \Prob_F\big(|\sfA(\mu_{F,\frakq})-\sfA(\overline{\nu})|>\Delta/2\big)
 \lesssim (M\Delta^2)^{-1}\lesssim_d M^{4/d-1} \lesssim (n \log n)^{4/d-1}.
\]
For $d>4$, the probability of the left-hand side is at most $0.05$ for all sufficiently large $n$.

Observe that $(\mathfrak{s}_b)_{\#}\overline{\nu}$ is supported in $\calX$. 
Now, let $\Lambda_{\frakq}$ be the distribution of $(\mu_{F,\frakq},(\mathfrak{s}_b)_{\#}\overline{\nu})$ on $\calP(\calX) \times \calP(\calX)$. Put 
\[
 c=(\underline{\delta}^2+3\Delta)^{1/2} \quad \text{and} \quad 
 s=\frac{(\underline{\delta}^2+4\Delta)^{1/2}-(\underline{\delta}^2+3\Delta)^{1/2}}{2}.
\]
On the event $|\sfA(\mu_{F,\frakq})-\sfA(\overline{\nu})|\le\Delta/2$,
\[
\GW^2(\mu_{F,\frakq},(\mathfrak{s}_b)_{\#}\overline{\nu}) \ge \underline{\delta}^2+\Delta/2. 
\]
For $\frakq\in\Dminus$, $\GW^2(\mu_{F,\frakq},(\mathfrak{s}_b)_{\#}\overline{\nu})$ is at most $\underline{\delta}^2+5\Delta/2$ whenever $\GW^2(\mu_{F,\frakq},\overline{\nu})\le\Delta$. For $\frakq\in\Dplus$, $\GW^2(\mu_{F,\frakq},(\mathfrak{s}_b)_{\#}\overline{\nu})$ is at least $\underline{\delta}^2+\Delta/2+16b^2\Delta\ge\underline{\delta}^2+4\Delta$. As such, from (\ref{eq: step 1}), we conclude that
\[
\begin{split}
&\Lambda_{\frakq}(\underline{\delta}\le\GW\le c)  \ge 0.85, \quad \frakq \in \Dminus, \\
&\Lambda_{\frakq}(\GW\ge c+2s) \ge 0.95, \quad \frakq \in \Dplus,
\end{split}
\]
where we used the fact that $\sup_{\frakq: \chi^2(\frakq,\fraku) \le 9} \Prob_F\big(|\sfA(\mu_{F,\frakq})-\sfA(\overline{\nu})|>\Delta/2\big) \le 0.05$. 

\medskip
\underline{Step 3} (Application of Proposition \ref{prop: fuzzy}).
Fix an estimator $\hat G$ and let $\psi=\ind_{\{\hat G\ge c+s\}}$.  As in the lower-bound proof for $1 \le p \le 2$ in Theorem \ref{thm: Wp main}, we may view $\psi$ as a randomized test based on  $I_1,\ldots,I_N\stackrel{\mathrm{i.i.d.}}{\sim}\frakq$.  Lemma \ref{lem: Wp testing} gives $\frakq^-\in\Dminus$ and $\frakq^+\in\Dplus$ such that
\[
 \Prob_{\frakq^-}(\psi=1)+\Prob_{\frakq^+}(\psi=0)\ge0.8.
\]

Since $\Lambda_{\frakq}$ need not be supported in $\{ \GW \ge \underline{\delta} \}$, we take 
 $\Lambda_0$ to be the conditional distribution of $\Lambda_{\frakq^-}$ given $\{\underline{\delta}\le\GW\le c\}$, and $\Lambda_1$ to be the conditional distribution of $\Lambda_{\frakq^+}$ given $\{\GW\ge c+2s\}$. Both priors are supported in $\{\GW\ge\underline{\delta}\}$, and Condition (i) of Proposition \ref{prop: fuzzy} holds with $\beta_0=\beta_1=0$. Let $\overline{P}_0,\overline{P}_1$ be the corresponding marginal distributions. Conditioning a prior on an event of probability $1-\varepsilon$ changes the marginal probability of any event by at most $\varepsilon$. As such, from Step 2, we conclude
\[
 \overline{P}_0(\hat G\ge c+s)+\overline{P}_1(\hat G<c+s)
 \ge \Prob_{\frakq^-}(\psi=1)+\Prob_{\frakq^+}(\psi=0) -0.15-0.05 \ge 0.6.
\]
Condition (ii) holds with $\alpha=0.4$. Proposition \ref{prop: fuzzy} now yields
\[
 \sup_{\GW(\mu,\nu)\ge\underline{\delta}}
 \Prob_{\mu,\nu}\big(|\hat G-\GW(\mu,\nu)|\ge s\big)\ge0.3.
\]
The conclusion follows as $s \asymp \Delta^{1/2}\wedge(\underline{\delta}^{-1}\Delta)
 \asymp_d (n\log n)^{-1/d}\wedge\big(\underline{\delta}^{-1}(n\log n)^{-2/d}\big)$. 
\end{proof}

\subsection{Upper bound}
Pick any $(\mu,\nu) \in \calP(B_{d_x}(0,1)) \times \calP(B_{d_y}(0,1))$ such that $t:=\GW(\mu,\nu) \ge \underline{\delta}$. 
Let $\hat G=\GW(\hat\mu_n,\hat\nu_m)$. By Theorem 4.2 in \cite{zhang2024gromov}, we have
\[
 \E[|\hat G^2-t^2|]\lesssim_{d_x,d_y}\overline{\varphi}_N^2.
\]
Since $\hat G,t\ge0$, we have $|\hat G-t|^2\le|\hat G^2-t^2|$. Also, for $t\ge\underline{\delta}>0$,
\[
 |\hat G-t| \le \underline{\delta}^{-1}|\hat G^2-t^2|.
\]
Taking expectation and applying Jensen's inequality, we obtain
\[
 \E[|\hat G-t|]\lesssim_{d_x,d_y}
 \overline{\varphi}_N\wedge(\underline{\delta}^{-1}\overline{\varphi}_N^2),
\]
which proves the upper bound in (\ref{eq: GW-main}). \qed

\appendix 

\section{On Remark \ref{rem: nonemptiness}}

This appendix concerns Remark \ref{rem: nonemptiness}. Consider the setting of Section \ref{sec: preliminaries}. For ease of notation, let $\calP_d^{ac}$ denote the set of all absolutely continuous $\mu \in \calP(B_d(0,1))$. 

\begin{lemma}
\label{lem: nonemptiness}
The following holds. 
\begin{enumerate}
\item[(i)] For any $p \ge 1$, $\{ \Wp (\mu,\nu) : (\mu,\nu) \in \calP_d^{ac} \times \calP_d^{ac} \} = [0,2)$.
\item[(ii)]
\[
\{ \GW (\mu,\nu) : (\mu,\nu) \in \calP_{d_x}^{ac} \times \calP_{d_y}^{ac} \} = 
\begin{cases}
[0,2\sqrt{2}), & d_x = d_y, \\
(0,2\sqrt{2}), & d_x \ne d_y. 
\end{cases}
\]
\end{enumerate}
\end{lemma}

\begin{proof}
(i). Let $v \in \R^d$ be a unit vector. For $0 \le t < 1$, set $r=1-t$, and let $\mu_t$ and $\nu_t$ be the uniform distributions on $B_d(tv,r)$ and $B_d(-tv,r)$, respectively. Translation by $2tv$ gives a coupling with cost $(2t)^p$. On the other hand, for any coupling $(X,Y)$, 
\[
\E[\|X-Y\|^p] \ge \| \E[X-Y] \|^p = (2t)^p,
\]
by Jensen's inequality.
Conclude that $\Wp(\mu_t,\nu_t) = 2t$.

Any absolutely continuous marginals charge zero mass on the boundary of $B_d(0,1)$, so any product coupling $(X,Y)$ satisfies $\|X-Y\| \le \|X\| + \|Y\| < 2$ almost surely, implying that $\Wp(\mu,\nu) < 2$ for any $\mu,\nu \in \calP_d^{ac}$.

\medskip

(ii). Let $0 < r < 1/3$, and let $v \in \R^{d_y}$ be a unit vector. Let $\mu_r$ be the uniform distribution on $B_{d_x}(0,r)$, and $\nu_r$ be the equal mixture of the uniform distributions on $B_{d_y}((1-r)v,r)$ and $B_{d_y}(-(1-r)v,r)$. 

From any coupling, generate independent draws $(X,Y), (X',Y')$. We have $\|X-X'\| \le 2r$, and that with probability $1/2$, $Y$ and $Y'$ fall in the opposite balls, giving $\|Y-Y'\|\ge2-4r>2r$. As such, using the representation (\ref{eq: GW-fourth-moments}), we obtain
\[
\GW(\mu_r,\nu_r)
\ge \frac{(2-4r)^2-(2r)^2}{\sqrt2},
\]
where the right-hand side is strictly decreasing between $0 < r < 1/3$, approaching $0$ as $r \to 1/3$, and approaching $2\sqrt{2}$ as $r \to 0$. To see that the range includes $(0,2\sqrt{2})$, pick any $0 < t < 2\sqrt{2}$, and choose $0 < r < 1/3$ so that $s := \GW(\mu_r,\nu_r) > t$. Dilate both marginals by $c:=\sqrt{t/s} < 1$. The dilated marginals are absolutely continuous and supported in the unit balls. Lemma \ref{lem: GW-scaling} yields $\GW((\mathfrak{s}_{c})_{\#}\mu_r,(\mathfrak{s}_{c})_{\#}\nu_r) = c^2s = t$.

To see that $2\sqrt{2}$ is not in the range, generate $(X,Y), (X',Y')$ from the independent coupling, and set $U=\|X-X'\|^2, V=\|Y-Y'\|^2, a = \E[U]$, and $b=\E[V]$. Since absolute continuity charges zero mass to the unit spheres, we have $a \le 2\E[\|X\|^2]<2$ and $b <2$. Since \(U,V\in[0,4]\) are independent,
\[
\GW^2(\mu,\nu) \le \E[(U-V)^2] \le 4a+4b-2ab =8-2(2-a)(2-b)<8.
\]

Finally, if $d_x = d_y$, then one may take $\mu = \nu$ and $\GW(\mu,\nu) = 0$. On the other hand, if $d_x \ne d_y$, then absolute continuity implies that the Hausdorff dimensions of $\spt(\mu)$ and $\spt(\nu)$ are $d_x$ and $d_y$, respectively.\footnote{See \cite{evans1992measure} for Hausdorff dimensions.} 
If $\GW(\mu,\nu)=0$,  then there would be a measure-preserving isometry between the supports. Isometries preserve Hausdorff dimension, which would force $d_x = d_y$, a contradiction. 

\end{proof}

\section*{Acknowledgments}

\textbf{Statement of AI use.} The authors used ChatGPT 5.6 Sol, 6.0 and 6.2 Astra (``AI models'') to assist with technical proofs, generating figures and tables, and proofreading. With instructions on proof strategies from the authors,  the AI models constructed first rough proofs for the main theorems, and the authors thoroughly revised and restructured the proofs. The writing is due to the authors. The authors take full responsibility for all aspects of the paper.

\bibliographystyle{alpha}
\bibliography{reference}

\end{document}